\documentclass{article}

\usepackage{arxiv}
\usepackage[utf8]{inputenc}
\usepackage[T1]{fontenc}
\usepackage{amsmath,amssymb,amsfonts,amsthm,mathtools}
\usepackage{bm}
\usepackage{graphicx}
\graphicspath{{./images/}}
\usepackage[dvipsnames]{xcolor}
\usepackage[numbers]{natbib}
\usepackage{hyperref}
\hypersetup{hidelinks,pdftitle={Resolvent-Based Singular-Value Diagnostics for Data-Driven Koopman Finite Sections}}
\usepackage{url}
\usepackage{nicefrac}
\usepackage{microtype}
\usepackage{multirow,array,booktabs,caption,subcaption}
\usepackage{enumitem}
\usepackage{algorithm}
\usepackage{algpseudocode}
\usepackage{float}
\usepackage{geometry}
\usepackage{afterpage}

\theoremstyle{plain}
\newtheorem{theorem}{Theorem}[section]

\newtheorem{proposition}[theorem]{Proposition}
\newtheorem{corollary}[theorem]{Corollary}
\newtheorem{lemma}[theorem]{Lemma}
\theoremstyle{definition}
\newtheorem{assumption}[theorem]{Assumption}
\newtheorem{definition}[theorem]{Definition}
\newtheorem{example}[theorem]{Example}
\theoremstyle{remark}
\newtheorem{remark}[theorem]{Remark}

\title{Resolvent-Based Singular-Value Diagnostics for Data-Driven Koopman Finite Sections}
\author{Yuanchao Xu$^{1}$\thanks{Corresponding author, xu.yuanchao.3a@kyoto-u.ac.jp}\enspace\thanks{These authors contributed equally.}\enspace \quad 
Itsushi Sakata$^{2}$\footnotemark[2]\enspace \quad 
Isao Ishikawa$^{1}$\enspace \\ 
$^1$ Center for Science Adventure and Collaborative Research Advancement (SACRA)\\
Graduate School of Science,
Kyoto University\\
$^2$ RIKEN Center for Advanced Intelligence Project\\
}

\begin{document}

\maketitle

\begin{abstract}
Finite-dimensional Koopman eigenvalues do not characterize resolvent growth, particularly for nonnormal compressions. We study the singular-value structure of empirical Koopman finite sections in the inner product induced by the data. Resolvent Dynamic Mode Decomposition (Resolvent DMD) removes unresolved Gram directions and then forms the Koopman or generator compression in a Gram-orthonormal basis. At a prescribed spectral parameter, the reciprocal of the smallest shifted singular value is the resolvent norm on the retained observable space. The associated right and left singular vectors give a minimum-residual pseudomode and the corresponding optimal forcing direction. For deterministic unregularized data, the squared ResDMD residual admits an orthogonal decomposition into a projected shifted residual and an invariance defect. The finite-section singular value is therefore a lower bound for the minimum ResDMD residual on the same space. We then identify stability assumptions that connect the finite sections with the underlying Koopman operator. Pointwise stability gives one-sided inclusion of inverse-resolvent sublevel sets. Stability on a punctured isolating neighborhood obtains local Hausdorff convergence near an isolated eigenvalue, while multiplicity stability preserves the algebraic count. Contractive, measure-preserving, and compact settings provide concrete sufficient conditions. The numerical examples examine coordinate dependence, finite-data residuals, contour calculations, and pseudomode separation for a noisy two-frequency signal.
\end{abstract}

%%%%%%%% Section 1 %%%%%%%%
\section{Introduction}
Many analyses of nonlinear dynamics begin with a differential equation or map and a numerical discretization \cite{strogatz2024nonlinear}. In applications where the governing equations are unavailable or uncertain, operator approximations must instead be constructed from observations. The resulting finite-dimensional operators act on a prescribed space of observables and depend on both the data and the approximation space.

The Koopman operator advances observables under the dynamics and remains linear when the state evolution is nonlinear \cite{Koopman1931,koopman1932dynamical,mezic2005spectral}. Dynamic Mode Decomposition (DMD) \cite{schmid2010dynamic,tu2013dynamic} and Extended Dynamic Mode Decomposition (EDMD) \cite{Williams_2015} replace this operator by finite matrices. Their convergence and spectral interpretation depend on the sampling rule, dictionary, observable space, and reference measure \cite{colbrook2024rigorous,ding2023higher,KLUS2020132416,korda2018convergence,mezic2020spectrum}. The eigenvalues of an EDMD matrix describe only part of its spectral behavior. The limiting Koopman operator may have continuous spectrum, while a nonnormal compression may have a large resolvent far from its eigenvalues. Computed eigenvalues may also vary with the dictionary, sampling rule, and regularization. Residual DMD (ResDMD) tests candidate spectral values through additional Gram matrices and thereby controls spectral pollution \cite{colbrook2023residual,colbrook2024rigorous,xu2025reskoopnet}. It does not, however, determine the left and right singular directions of the shifted compression.

We consider the resolvent of the empirical Koopman compression in Gram-orthonormal coordinates. For each prescribed \(z\in\mathbb C\), the smallest singular triplet of the shifted Koopman or generator matrix yields the resolvent norm, a minimum-residual pseudomode, and its optimal forcing direction. We use the term Resolvent DMD for this family of shifted singular-value calculations. Unlike input--output resolvent analysis for linearized state equations \cite{herrmann2021data}, the present resolvent acts on observables and is evaluated in the complex spectral plane.

Gram coordinates are needed when the dictionary is nonorthogonal or numerically rank deficient. Rank truncation is imposed before the shifted matrices are formed, so unresolved directions do not enter the pseudoinverse. On a common retained subspace, the deterministic residual splits orthogonally into a projected shifted residual and an invariance defect. The finite-section singular value consequently supplies a necessary condition for a small residual.

The passage from finite sections to the Koopman operator is governed by resolvent stability. Pointwise stability implies retained strong resolvent convergence and one-sided inclusion of inverse-resolvent sublevel sets. Stability on a punctured isolating region gives local Hausdorff convergence near an isolated eigenvalue. A separate multiplicity condition recovers the algebraic count. The backward-shift example shows that resolvent sublevel sets may retain spectral information that is absent from the finite-section eigenvalues. The numerical examples treat coefficient sensitivity, finite-data residuals, contour fields, and pseudomode-based separation. Algorithm~\ref{alg:resolvent_dmd} gives the corresponding Gram-coordinate computation.

Section~\ref{sec:background} introduces the Koopman and resolvent quantities used in the analysis. Section~\ref{sec:method} develops the Gram-weighted finite-section construction and its relation to ResDMD. Section~\ref{sec:cvg_analysis} treats finite-section stability, and Section~\ref{sec:numerical_studies} presents four numerical examples. Sections~\ref{sec:discussion} and~\ref{sec:conclusion} discuss the scope of the results and summarize the main conclusions. The proofs and supplementary pendulum calculations are given in the appendices.

%%%%%%%% Section 2 %%%%%%%%
\section{Koopman operators and resolvent analysis}\label{sec:background}
%%%%%%%% Section 2.1 %%%%%%%%
\subsection{Koopman operator}
Let \(\mathcal M\) be a state space with measure \(\mu\), and let \(T\colon\mathcal M\to\mathcal M\) be a discrete-time map. Assume that \(T_*\mu\ll\mu\) and that \(d(T_*\mu)/d\mu\) is essentially bounded. Composition with \(T\) then defines a bounded operator on \(\mathcal F=L^2(\mathcal M,\mu)\). With
\[
    \langle f,g\rangle=\int_{\mathcal M}\overline f g\,d\mu,
\]
the Koopman operator is
\[
    \mathcal K f(x)=f(T(x)),\qquad f\in\mathcal F.
\]

For continuous-time dynamics, let \(\{T_t\}_{t\ge0}\) be a flow and define
\[
    \mathcal K^t f(x)=f(T_t(x)),\qquad t\ge0.
\]
Under the standard \(C_0\)-semigroup assumptions \cite{engel1999one,pazy2012semigroups}, \(\{\mathcal K^t\}_{t\ge0}\) has a closed, densely defined generator \(\mathcal A\) and \(\mathcal K^t=e^{t\mathcal A}\). For a fixed sampling interval \(\tau>0\), write
\begin{equation}\label{eq:embedding}
    \mathcal K=\mathcal K^\tau.
\end{equation}
In the numerical examples, \(\tau\) is the stated sampling interval \(\Delta t\). We treat the sampled operator and the generator as distinct objects. Spectral mapping is used only when its hypotheses are satisfied. If \(\mathcal K\) is contractive, then \(\sigma(\mathcal K)\) lies in the closed unit disk. If the Koopman semigroup is contractive on the chosen observable space, then \(\sigma(\mathcal A)\) lies in the closed left half-plane. Dissipation of the state dynamics alone does not imply either operator statement. The norm also depends on the observable space and the reference measure.

For a general \(C_0\)-semigroup, only the inclusion
\[
    e^{t\sigma(\mathcal A)}\subseteq\sigma(\mathcal K^t)\setminus\{0\}
\]
is guaranteed. Equality may fail without additional assumptions \cite{engel1999one,pazy2012semigroups}. We therefore evaluate the generator resolvent directly rather than infer it from the sampled Koopman spectrum. The sampled operator can be approximated from snapshot pairs without an explicit formula for \(T\) \cite{budivsic2012applied,KLUS2020132416,Williams_2015}.

%%%%%%%% Section 2.2 %%%%%%%%
\subsection{Resolvents, shifted singular values, and pseudomodes}
Let \(\mathcal L\) be a closed, densely defined operator on \(\mathcal F\). Its resolvent set is \(\rho(\mathcal L)=\mathbb C\setminus\sigma(\mathcal L)\), and
\[
    \mathcal R(z,\mathcal L)=(zI-\mathcal L)^{-1},
    \qquad z\in\rho(\mathcal L).
\]
The resolvent is analytic on \(\rho(\mathcal L)\) and satisfies
\[
    \mathcal R(z,\mathcal L)-\mathcal R(w,\mathcal L)
    =(w-z)\mathcal R(z,\mathcal L)\mathcal R(w,\mathcal L).
\]
For \(\varepsilon>0\), the \(\varepsilon\)-pseudospectrum is
\[
\begin{aligned}
    \sigma_\varepsilon(\mathcal L)
    &=\sigma(\mathcal L)\cup
      \bigl\{z\in\rho(\mathcal L)\mid
      \|\mathcal R(z,\mathcal L)\|>\varepsilon^{-1}\bigr\}.
\end{aligned}
\]
For a finite matrix \(L\), this is equivalently
\[
    \sigma_\varepsilon(L)
    =\bigl\{z\in\mathbb C\mid
      \sigma_{\min}(zI-L)<\varepsilon\bigr\}.
\]
A grid calculation of the resolvent therefore reduces to independent smallest-singular-value problems for the shifted matrices \cite{trefethen2005spectra}. Let \((s,u,v)\) be a smallest singular triplet of \(zI-L\) satisfying
\[
    (zI-L)v=su,\qquad
    \|u\|_2=\|v\|_2=1,\qquad
    s=\sigma_{\min}(zI-L).
\]
The right singular vector \(v\) minimizes \(\|(zI-L)w\|_2\) over unit vectors \(w\) and is therefore a minimum-residual pseudomode. If \(s>0\), then
\[
    (zI-L)^{-1}u=s^{-1}v,
    \qquad
    \|(zI-L)^{-1}\|_2=s^{-1}.
\]
Thus \(u\) is an optimal forcing direction and \(v\) is the corresponding response. For a normal matrix, \(\sigma_{\min}(zI-L)=\operatorname{dist}(z,\sigma(L))\). For a nonnormal matrix, the singular value may be much smaller than the spectral distance. Large resolvent values can then occur away from the eigenvalues. If \(\mathcal K\) is contractive, the Neumann series gives \(\|\mathcal R(z,\mathcal K)\|\le (|z|-1)^{-1}\) for \(|z|>1\). We write \(\Gamma_R=\{z\in\mathbb C\mid |z|=R\}\) for the circle of radius \(R>0\).

%%%%%%%% Section 3 %%%%%%%%
\section{Gram-weighted resolvent analysis for Koopman finite sections}\label{sec:method}
Let \(\psi_1,\ldots,\psi_N\) be a dictionary and let \((x_m,y_m)_{m=1}^M\) be snapshot pairs. The matrices \(\Psi_X,\Psi_Y\in\mathbb C^{M\times N}\) contain the dictionary values at the initial and successor states. Let
\[
    W=\operatorname{diag}(w_1,\ldots,w_M),\qquad
    w_m\ge0,\qquad \sum_{m=1}^M w_m=1,
\]
be a sampling or quadrature weight matrix. Uniform empirical weights correspond to \(W=M^{-1}I_M\). Define
\begin{equation}\label{eq:edmd}
\begin{aligned}
    G_{N,M}&=\Psi_X^*W\Psi_X,&
    B_{N,M}&=\Psi_X^*W\Psi_Y,
    &\mathbf K_{N,M}&=G_{N,M}^\dagger B_{N,M},\\
    B^{\mathcal A}_{N,M}&=\Psi_X^*W\Psi_X',&&
    &\mathbf A_{N,M}&=G_{N,M}^\dagger B^{\mathcal A}_{N,M}.
\end{aligned}
\end{equation}
Here \(\Psi_X'\) contains the generator action on the dictionary. When exact generator data are unavailable, one may use \(\Psi_X'=(\Psi_Y-\Psi_X)/\Delta t\). For stochastic dynamics, the generator action may instead be evaluated with It\^o's formula \cite{KLUS2020132416,xu2025datadrivenframeworkkoopmansemigroup,xu2026reinforcement}. We reserve \(\mathbf K_N\) and \(\mathbf A_N\) for fixed-\(N\) population matrices. The subscript \((N,M)\) denotes a finite-data quantity.

\begin{remark}\label{rm:fd_accuracy}
The approximation \(\Psi_X'=(\Psi_Y-\Psi_X)/\Delta t\) is first-order accurate. If \(\mathcal A\varphi=\lambda\varphi\), then
\[
\frac{\mathcal K^{\Delta t}-I}{\Delta t}\varphi
=\left(\lambda+\frac{\lambda^2\Delta t}{2}+O(\Delta t^2)\right)\varphi.
\]
For \(\lambda=i\omega\), the leading artificial real part is \(-\omega^2\Delta t/2\). The relative frequency error is \(O((\omega\Delta t)^2)\), with leading coefficient \(1/6\). Exact derivatives \(\nabla\psi_j\cdot F\) avoid this temporal discretization error when they are available \cite{KLUS2020132416}. Decreasing \(\Delta t\) reduces truncation error but may amplify observational noise in the finite difference.
\end{remark}

The resolvent in coefficient coordinates is
\begin{equation}\label{eq:resolvent_matrix_deterministic}
    \mathbf R_{N,M}(z)=(zI-\mathbf K_{N,M})^{-1},
    \qquad z\in\rho(\mathbf K_{N,M}).
\end{equation}
The Euclidean coefficient norm agrees with the empirical observable norm only in empirically orthonormal coordinates. Choose a numerical rank \(r\) and write
\[
    G_{N,M}=V_r\Sigma_rV_r^*+G_{\rm disc},
    \qquad \Sigma_r=\operatorname{diag}(\sigma_1,\ldots,\sigma_r)>0.
\]
The columns of \(V_r\) are orthonormal eigenvectors associated with the retained eigenvalues of \(G_{N,M}\). They satisfy \(V_r^*G_{\rm disc}V_r=0\). The coordinate transformation
\[
    c=V_r\Sigma_r^{-1/2}q
\]
identifies Euclidean norm in \(q\) with the empirical observable norm. The discrete-time and generator compressions in these coordinates are
\[
    \widetilde K_{N,M}
    =\Sigma_r^{-1/2}V_r^*B_{N,M}V_r\Sigma_r^{-1/2},
    \qquad
    \widetilde A_{N,M}
    =\Sigma_r^{-1/2}V_r^*B^{\mathcal A}_{N,M}V_r\Sigma_r^{-1/2}.
\]
If \(G_{N,M}\) is nonsingular, one may take \(r=N\). This construction is a rank-truncated compression. Replacing the pseudoinverse by a Tikhonov filter gives a different regularized operator.

For the discrete-time problem, define
\begin{equation}\label{eq:gram_resolvent}
    \widetilde S_{N,M}(z)=zI_r-\widetilde K_{N,M},
    \qquad
    s_{N,M}(z)=\sigma_{\min}\bigl(\widetilde S_{N,M}(z)\bigr).
\end{equation}
The generator analogue is
\begin{equation}\label{eq:generator_resolvent}
    \widetilde R^{\mathcal A}_{N,M}(z)
    =(zI_r-\widetilde A_{N,M})^{-1}.
\end{equation}

\begin{proposition}[Gram-geometry identity]\label{prop:gram_identity}
The columns of \(\Phi_{X,r}=\Psi_XV_r\Sigma_r^{-1/2}\) are orthonormal in the \(W\)-inner product. In this basis, \(\widetilde K_{N,M}\) represents the retained empirical finite section. For \(z\in\rho(\widetilde K_{N,M})\),
\begin{equation}\label{eq:gram_identity}
    \bigl\|(zI_r-\widetilde K_{N,M})^{-1}\bigr\|_2
    =s_{N,M}(z)^{-1}.
\end{equation}
If \(G_{N,M}\) is nonsingular, then
\[
    s_{N,M}(z)=\sigma_{\min}\!\left(
      G_{N,M}^{1/2}(zI-\mathbf K_{N,M})G_{N,M}^{-1/2}
    \right).
\]
Moreover, \(s_{N,M}(z)\) is invariant under unitary changes of basis within the retained empirical subspace.
\end{proposition}

Thus
\[
    \{z\in\mathbb C\mid s_{N,M}(z)<\varepsilon\}
\]
is the \(\varepsilon\)-pseudospectrum of \(\widetilde K_{N,M}\). Algorithm~\ref{alg:resolvent_dmd} evaluates this set only at the prescribed probe points.

\begin{remark}\label{rm:coordinate_conditioning}
Suppose that \(G_{N,M}\) is nonsingular and set \(\kappa=\kappa_2(G_{N,M})\). Then
\[
  \kappa^{-1/2}s_{N,M}(z)
  \le \sigma_{\min}(zI-\mathbf K_{N,M})
  \le \kappa^{1/2}s_{N,M}(z).
\]
The bounds follow from the similarity relation in Proposition~\ref{prop:gram_identity}. Ill-conditioned dictionary coordinates may therefore distort coefficient-space singular values and the associated absolute thresholds.
\end{remark}

Rank truncation precedes the assembly of the shifted matrices. Poorly resolved directions are therefore removed before they enter the pseudoinverse.

Let
\[
    \widetilde S_{N,M}(z)v_z=s_{N,M}(z)u_z,
    \qquad \|u_z\|_2=\|v_z\|_2=1.
\]
Then \(v_z\) is a minimum-residual pseudomode in retained coordinates. The corresponding observable has dictionary coefficients
\[
    c_z=V_r\Sigma_r^{-1/2}v_z,
    \qquad
    g_z=\sum_j(c_z)_j\psi_j.
\]
If \(s_{N,M}(z)>0\), then
\[
    (zI_r-\widetilde K_{N,M})^{-1}u_z
    =s_{N,M}(z)^{-1}v_z,
\]
so \(u_z\) is the associated optimal forcing direction.

ResDMD evaluates Koopman residuals with an additional output Gram matrix \cite{colbrook2023residual,colbrook2024rigorous}. The following proposition relates this residual to the shifted singular value on the retained space.

\begin{proposition}\label{prop:resdmd_gap}
Assume deterministic data without regularization. Use the same weights and retained subspace in the finite-section and ResDMD calculations. Let \(\|h\|_W^2=h^*Wh\), and let \(P_r\) be the \(W\)-orthogonal projection onto \(\operatorname{ran}(\Phi_{X,r})\). For \(c=V_r\Sigma_r^{-1/2}q\) with \(\|q\|_2=1\), or equivalently \(\|\Psi_Xc\|_W=1\),
\begin{equation}\label{eq:resdmd_decomposition}
\begin{aligned}
    \|\Psi_Yc-z\Psi_Xc\|_W^2
    &=\|(zI_r-\widetilde K_{N,M})q\|_2^2
      +\|(I-P_r)\Psi_Yc\|_W^2.
\end{aligned}
\end{equation}
If
\[
  \tau_{N,M}(z)
  \coloneqq
  \min_{\|q\|_2=1}
  \|\Psi_YV_r\Sigma_r^{-1/2}q
     -z\Psi_XV_r\Sigma_r^{-1/2}q\|_W,
\]
then
\[
  \tau_{N,M}(z)\ge s_{N,M}(z),
  \qquad
  \{z\mid\tau_{N,M}(z)<\varepsilon\}
  \subseteq
  \{z\mid s_{N,M}(z)<\varepsilon\}.
\]
\end{proposition}

The second term in Eq.~\eqref{eq:resdmd_decomposition} is the empirical invariance defect in the direction \(c\). Hence \(s_{N,M}(z)\ge\varepsilon\) rules out a ResDMD residual below \(\varepsilon\) on the same retained subspace. The singular value gives a necessary test, while the defect measures the component of the propagated observable outside the trial space. Equation~\eqref{eq:resdmd_decomposition} changes under regularization and for stochastic data with one successor per state. In either case, the corresponding inequality must be proved for the estimator under consideration. Use \(\Psi_Z=\Psi_Y\) in discrete time and \(\Psi_Z=\Psi_X'\) for the generator calculation. Algorithm~\ref{alg:resolvent_dmd} states the computation in retained coordinates.

\begin{algorithm}[t]
\caption{Evaluation of retained finite-section resolvent quantities in Gram coordinates}
\label{alg:resolvent_dmd}
\begin{algorithmic}[1]
\Require Snapshot matrix \(\Psi_X\), companion matrix \(\Psi_Z=\Psi_Y\) for discrete time or \(\Psi_Z=\Psi_X'\) for a generator, optional weights \(W\), spectral points \(\{z_k\}_{k=1}^K\), threshold \(\varepsilon>0\), retained-rank rule, and numerical tolerance \(\tau_{\rm num}\).
\State If \(W\) is omitted, set \(W=M^{-1}I_M\). Form \(G=\Psi_X^*W\Psi_X\) and \(B=\Psi_X^*W\Psi_Z\).
\State Compute \(G=V_r\Sigma_rV_r^*+G_{\rm disc}\) according to the retained-rank rule. Set \(r=N\) if no directions are discarded.
\State Form \(\widetilde L=\Sigma_r^{-1/2}V_r^*BV_r\Sigma_r^{-1/2}\).
\For{\(k=1,\ldots,K\)}
  \State Set \(S_k=z_kI_r-\widetilde L\) and compute a smallest singular triplet \(S_kv_k=s_ku_k\), without forming \(S_k^{-1}\).
  \State If \(s_k<\varepsilon\), retain \(z_k\) in the inverse-resolvent sublevel set. If \(s_k\le\tau_{\rm num}\), mark the value as numerically unresolved.
  \State Set \(c_k=V_r\Sigma_r^{-1/2}v_k\) and \(g_k=\sum_j(c_k)_j\psi_j\).
\EndFor
\State \Return \(\{z_k\mid s_k<\varepsilon\}\), \(\{s_k\}_{k=1}^K\), pseudomodes \(\{g_k\}_{k=1}^K\), and numerical-resolution flags.
\end{algorithmic}
\end{algorithm}

The returned observables \(g_k\) are the pseudomodes expressed in the original dictionary. Forming \(G\) and \(B\) requires \(O(MN^2)\) operations, while a dense eigendecomposition of \(G\) requires \(O(N^3)\). A dense smallest-singular-triplet calculation costs \(O(r^3)\) at each spectral point. Iterative solvers may be used when only the smallest triplet is required. The shifted problems are independent and may be solved in parallel.

The calculation is exact for the retained empirical finite section. Section~\ref{sec:cvg_analysis} gives conditions that relate these quantities to the resolvent and spectrum of the limiting operator. The discrete-time and generator problems remain separate throughout the analysis.

%%%%%%%% Section 4 %%%%%%%%
\section{Operator interpretation under finite-section stability}\label{sec:cvg_analysis}
The identities in Section~\ref{sec:method} hold on the retained empirical subspace without approximation. Their relation to the spectrum and resolvent of the Koopman operator requires stability of the finite sections. We state the assumptions and convergence results below. The proofs are given in Appendix~\ref{app:proofs}.

\begin{definition}[Finite-section setting]\label{def:finite_section_resolvents}
Let \(\mathcal K\) be a bounded Koopman operator on \(\mathcal F\). Let \(\mathcal F_N\subset\mathcal F\) be nested finite-dimensional Galerkin subspaces with orthogonal projections \(\Pi_N\) such that \(\Pi_N\to I\) strongly. The ideal finite section is
\[
    \mathcal K_N\coloneqq\Pi_N\mathcal K\Pi_N.
\]
The same symbol denotes the finite-rank operator on \(\mathcal F\) and its restriction to \(\mathcal F_N\). We write \(\mathbf K_N\) for the population matrix and \(\mathbf K_{N,M}\) for the empirical EDMD matrix. In this section, \(\sigma(\mathcal K_N)\) denotes the spectrum of \(\mathcal K_N|_{\mathcal F_N}\).

The limiting resolvent is
\[
    \mathcal R(z)\coloneqq\mathcal R(z,\mathcal K)
    =(zI-\mathcal K)^{-1},
    \qquad z\in\rho(\mathcal K).
\]
Whenever \(z\in\rho(\mathcal K_N|_{\mathcal F_N})\), define the retained inverse by
\[
    \mathcal R_N^{\rm ret}(z)
    \coloneqq(zI_{\mathcal F_N}-\mathcal K_N|_{\mathcal F_N})^{-1}.
\]
Its matrix representation is the large-data limit of Eq.~\eqref{eq:resolvent_matrix_deterministic}. We also write \(\mathcal R_N(z)\) for this operator on \(\mathcal F_N\). The projection \(\Pi_N\) is shown explicitly when the inverse acts on a vector in \(\mathcal F\).
\end{definition}

For fixed \(N\), convergence of \(\mathbf K_{N,M}\) to \(\mathbf K_N\) also requires stability of the retained rank. Appendix~\ref{app:finite_data_limit} summarizes the finite-data assumptions. The results below concern the Galerkin limit \(N\to\infty\). All finite-section resolvent norms are taken on \(\mathcal F_N\), or equivalently in the Gram geometry of Section~\ref{sec:method}.

\begin{assumption}[Pointwise finite-section stability]\label{ass:stable_cvg}
Let \(z\in\rho(\mathcal K)\). The finite sections are stable at \(z\) if there exist \(N(z)\in\mathbb N\) and \(M_z<\infty\) such that, for every \(N>N(z)\), the point \(z\) belongs to \(\rho(\mathcal K_N|_{\mathcal F_N})\) and
\begin{equation}\label{eq:stable_cvg}
    \|\mathcal R_N(z)\|\le M_z.
\end{equation}
Stability on a set \(U\subset\rho(\mathcal K)\) means stability at every point of \(U\).
\end{assumption}

Assumption~\ref{ass:stable_cvg} is the standard local stability condition for finite-section spectral approximation \cite[Chapter 3.4]{chaitin1983spectral}. It excludes unbounded growth of the finite-section resolvent at the specified points of \(\rho(\mathcal K)\).

Contour integrals require stability along the contour. Theorem~\ref{thm:hausdorff_cvg_isolated_eig} further requires stability throughout the punctured isolating region. In orthonormal coordinates, pointwise stability is equivalent to a uniform positive lower bound on \(\sigma_{\min}(zI-\mathbf K_N)\). For a nonorthogonal dictionary, the corresponding quantity is the Gram-weighted singular value from Section~\ref{sec:method}.

\begin{remark}\label{rm:retained_resolvent}
For the ideal Galerkin compressions, Assumption~\ref{ass:stable_cvg} implies retained strong resolvent convergence at every stable point.
\begin{equation}\label{eq:resolvent_strong_cvg}
    \mathcal R_N^{\rm ret}(z)\Pi_N f
    \longrightarrow\mathcal R(z)f,
    \qquad f\in\mathcal F.
\end{equation}
The derivation is given in Appendix~\ref{app:finite_data_limit}. For more general finite-section approximations, Eq.~\eqref{eq:resolvent_strong_cvg} must be assumed or verified independently.
\end{remark}

\begin{definition}\label{def:spec_proj}
Let \(\lambda\) be an isolated eigenvalue of \(\mathcal K\) with finite algebraic multiplicity \(m_\lambda\). Let \(\Gamma_\lambda\subset\rho(\mathcal K)\) be a positively oriented rectifiable Jordan curve whose interior contains \(\lambda\) and no other point of \(\sigma(\mathcal K)\). The associated Riesz projections \cite{kato1995perturbation} are
\[
    \mathcal P^\lambda
    \coloneqq\frac{1}{2\pi i}\int_{\Gamma_\lambda}\mathcal R(z)\,dz,
    \qquad
    \mathcal P_N^\lambda
    \coloneqq\frac{1}{2\pi i}\int_{\Gamma_\lambda}
    \mathcal R_N^{\rm ret}(z)\Pi_N\,dz,
\]
provided that \(\Gamma_\lambda\subset\rho(\mathcal K_N|_{\mathcal F_N})\). Their ranges are the corresponding invariant spectral subspaces. On \(\mathcal F_N\), the second integral is the usual Riesz projection. The factor \(\Pi_N\) extends it to \(\mathcal F\).
\end{definition}

\begin{assumption}[Multiplicity stability relative to an isolating contour]\label{ass:strong_stable_cvg}
Fix \(\lambda\) and \(\Gamma_\lambda\) as in Definition~\ref{def:spec_proj}. Suppose that Assumption~\ref{ass:stable_cvg} holds on \(\Gamma_\lambda\), so \(\mathcal P_N^\lambda\) is defined for all sufficiently large \(N\). We assume that
\[
    \dim\mathcal P_N^\lambda\mathcal F
    =\dim\mathcal P^\lambda\mathcal F
    =m_\lambda
\]
for all sufficiently large \(N\).
\end{assumption}

This is the multiplicity condition used in strongly stable spectral approximation \cite[Chapter 5]{chaitin1983spectral}. Corollary~\ref{cor:spectral_projection} gives the lower bound \(\dim\mathcal P_N^\lambda\mathcal F\ge m_\lambda\). Equality excludes additional finite-section eigenvalues inside the contour. Proposition~\ref{prop:compact_case} derives equality from norm convergence for compact operators.

\begin{lemma}\label{lem:uniform_bound}
Let \(U\Subset\rho(\mathcal K)\) be compact. If Assumption~\ref{ass:stable_cvg} holds at every \(z\in U\), then there exist \(N_U\in\mathbb N\) and \(M_U<\infty\) such that, for every \(N>N_U\),
\[
    U\subset\rho(\mathcal K_N|_{\mathcal F_N}),
    \qquad
    \sup_{z\in U}\|\mathcal R_N(z)\|\le M_U.
\]
\end{lemma}

\begin{corollary}\label{cor:spectral_projection}
Let \(\lambda\) be an isolated eigenvalue of \(\mathcal K\) with algebraic multiplicity \(m_\lambda\), and let \(\Gamma_\lambda\) be the isolating contour from Definition~\ref{def:spec_proj}. If Assumption~\ref{ass:stable_cvg} holds on \(\Gamma_\lambda\), then
\[
    \dim\bigl(\mathcal P_N^\lambda\mathcal F\bigr)
    \ge\dim\bigl(\mathcal P^\lambda\mathcal F\bigr)
    =m_\lambda
\]
for all sufficiently large \(N\).
\end{corollary}

\begin{theorem}\label{thm:hausdorff_cvg_isolated_eig}
Let \(\lambda\) be an isolated eigenvalue of \(\mathcal K\) with finite algebraic multiplicity \(m\). Let \(\Gamma\subset\rho(\mathcal K)\) be a positively oriented rectifiable Jordan curve enclosing \(\lambda\), and let \(\Delta\) denote its interior. Suppose that
\[
    \sigma(\mathcal K)\cap\overline\Delta=\{\lambda\}
\]
and that Assumption~\ref{ass:stable_cvg} holds on every compact subset of \(\overline\Delta\setminus\{\lambda\}\). Then \(\sigma(\mathcal K_N)\cap\Delta\ne\varnothing\) for all sufficiently large \(N\), and
\[
    d_H\bigl(\sigma(\mathcal K_N)\cap\Delta,\{\lambda\}\bigr)
    \longrightarrow0.
\]
Here \(d_H\) is the Hausdorff distance between nonempty compact subsets of \(\mathbb C\). If Assumption~\ref{ass:strong_stable_cvg} also holds relative to \(\Gamma\), then \(\sigma(\mathcal K_N)\cap\Delta\) contains exactly \(m\) eigenvalues, counted with algebraic multiplicity, for all sufficiently large \(N\).
\end{theorem}

\subsection{One-sided inverse-resolvent inclusion and eigenvalue loss}\label{subsec:pseudo_detection_limit}
Theorem~\ref{thm:hausdorff_cvg_isolated_eig} is local to an isolated eigenvalue. For a nonisolated spectral set, the finite-section eigenvalues may fail to approximate the spectrum even under retained strong resolvent convergence. The inverse-resolvent sublevel sets still satisfy a one-sided inclusion.

\begin{theorem}[One-sided inverse-resolvent inclusion]\label{thm:one_sided_indicator}
Suppose that Assumption~\ref{ass:stable_cvg} holds at \(z\in\rho(\mathcal K)\). Then
\[
    \liminf_{N\to\infty}\|\mathcal R_N(z)\|
    \ge\|\mathcal R(z)\|.
\]
Consequently, if \(\varepsilon>0\) and \(\|\mathcal R(z)\|>\varepsilon^{-1}\), then
\[
    \|\mathcal R_N(z)\|^{-1}<\varepsilon
\]
for all sufficiently large \(N\). Every stable point with \(\|\mathcal R(z)\|^{-1}<\varepsilon\) therefore enters the finite-section inverse-resolvent sublevel set for all sufficiently large \(N\).

Let \(Q\Subset\rho(\mathcal K)\) be compact and suppose that Assumption~\ref{ass:stable_cvg} holds on \(Q\). If
\[
    \inf_{z\in Q}\|\mathcal R(z)\|>\varepsilon^{-1},
\]
then
\[
    \sup_{z\in Q}\|\mathcal R_N(z)\|^{-1}<\varepsilon
\]
for all sufficiently large \(N\).
\end{theorem}

\begin{example}[Eigenvalue loss for finite sections of a shift]\label{ex:backward_shift_limitation}
Let \(\mathcal M=\mathbb N\) with counting measure and \(T(n)=n+1\). Then
\[
    \mathcal Kf(n)=f(n+1)
\]
is the backward shift on \(\ell^2(\mathbb N)\), a contraction with
\[
    \sigma(\mathcal K)=\overline{\mathbb D}.
\]
For the standard basis, let \(\mathcal F_N=\operatorname{span}\{e_1,\ldots,e_N\}\). Although \(\Pi_N\mathcal K\Pi_N\to\mathcal K\) strongly, each finite-section matrix is nilpotent, and hence
\[
    \sigma(\mathcal K_N)=\{0\}
    \qquad\text{for every }N.
\]
For \(|z|>1\), Proposition~\ref{prop:contractive_case} gives
\[
    \|\mathcal R_N(z)\|\le(|z|-1)^{-1},
\]
so the finite sections are stable outside the unit disk. Their eigenvalues nevertheless do not approximate \(\overline{\mathbb D}\). Indeed, if \(\Gamma=\{|z|=2\}\), \(\Delta=\{|z|<2\}\), and \(U=\overline{\mathbb D}\), then
\[
    \sup_{u\in U}\operatorname{dist}
    \bigl(u,\sigma(\mathcal K_N)\cap\Delta\bigr)=1
\]
for every \(N\).

The inverse-resolvent diagnostic has different behavior. For \(z\ne0\) and the final basis vector \(e_N\),
\[
  \|\mathcal R_N(z)e_N\|^2
  =\sum_{k=1}^{N}|z|^{-2k}.
\]
Hence, for \(0<|z|<1\),
\[
  \|\mathcal R_N(z)\|^{-1}\le |z|^N\longrightarrow0,
\]
and for \(|z|=1\),
\[
  \|\mathcal R_N(z)\|^{-1}\le N^{-1/2}\longrightarrow0.
\]
The point \(z=0\) is an eigenvalue of every finite section. Hence every fixed point of \(\overline{\mathbb D}\) eventually lies in any positive inverse-resolvent sublevel set, although all finite-section eigenvalues remain at the origin. The conclusion is one-sided. A fixed positive threshold may also include points outside the spectrum.
\end{example}

\begin{remark}
The shift is a standard example of nonnormal spectral approximation. For these Toeplitz truncations, fixed-threshold pseudospectra converge even though the truncated eigenvalues do not \cite{reichel1992eigenvalues}.
\end{remark}

Theorem~\ref{thm:one_sided_indicator} uses only finite-section stability. We next give standard settings in which this stability follows from contractivity, skew-adjointness, or compactness.

\subsection{Standard cases of resolvent stability}\label{subsec:verifiable_conditions}
Let \(\Pi_N\) be the \(L^2(\mu)\)-orthogonal projection onto a nested Galerkin space \(\mathcal F_N\). The following propositions cover three standard cases. They give direct resolvent bounds or, in the compact case, norm convergence. Nonorthogonal dictionaries are interpreted through the Gram coordinates of Section~\ref{sec:method}.

\begin{proposition}[Contractive discrete-time compressions]\label{prop:contractive_case}
Let \(\mathcal K_N=\Pi_N\mathcal K\Pi_N\). If \(\|\mathcal K\|\le1\), then, for every \(|z|>1\) and every \(N\),
\[
    z\in\rho(\mathcal K_N|_{\mathcal F_N}),
    \qquad
    \|\mathcal R_N(z)\|\le\frac{1}{|z|-1}.
\]
For a Koopman composition operator, contractivity holds, for example, when the reference measure is invariant or sub-invariant under the sampled dynamics, \(T_*\mu\le\mu\). Hence Assumption~\ref{ass:stable_cvg} holds uniformly in \(N\) for every \(|z|>1\).
\end{proposition}

\begin{proposition}[Generators of measure-preserving flows]\label{prop:generator_measure_preserving}
Let \(\{T_t\}_{t\in\mathbb R}\) be an invertible \(\mu\)-preserving flow. Suppose that \(\{\mathcal K^t\}_{t\in\mathbb R}\) is a strongly continuous unitary group on \(L^2(\mu)\) with skew-adjoint generator \(\mathcal A\). Let \(\mathcal F_N\subset D(\mathcal A)\) be nested finite-dimensional spaces with orthogonal projections \(\Pi_N\). Assume that \(\bigcup_N\mathcal F_N\) is a core for \(\mathcal A\), and define \(\mathcal A_N=\Pi_N\mathcal A\Pi_N\) on \(\mathcal F_N\). For \(\operatorname{Re} z\ne0\), set
\[
    \mathcal R_N^{\mathcal A}(z)
    \coloneqq\mathcal R^{\rm ret}(z,\mathcal A_N)\Pi_N,
    \qquad
    \mathcal R^{\rm ret}(z,\mathcal A_N)
    \coloneqq(zI_{\mathcal F_N}-\mathcal A_N|_{\mathcal F_N})^{-1}.
\]
The following statements hold.
\begin{enumerate}[label=(\roman*)]
    \item The compression \(\mathcal A_N\) is skew-adjoint on \(\mathcal F_N\). Hence \(\sigma(\mathcal A_N)\subset i\mathbb R\), and
    \[
        \|\mathcal R^{\rm ret}(z,\mathcal A_N)\|
        =\frac{1}{\operatorname{dist}(z,\sigma(\mathcal A_N))}
        \le\frac{1}{|\operatorname{Re} z|}
    \]
    for every \(\operatorname{Re} z\ne0\), uniformly in \(N\).
    \item For every \(\operatorname{Re} z\ne0\),
    \[
        \mathcal R_N^{\mathcal A}(z)h
        \longrightarrow\mathcal R(z,\mathcal A)h,
        \qquad h\in L^2(\mu).
    \]
    \item For every \(\lambda\in\sigma(\mathcal A)\), including points in the continuous spectrum and the possible point \(0\),
    \[
        \operatorname{dist}(\lambda,\sigma(\mathcal A_N))\longrightarrow0.
    \]
\end{enumerate}
If \(\{\mathcal K^t\}_{t\ge0}\) is instead a contraction semigroup, then \(\mathcal A_N\) is dissipative and
\[
    \|\mathcal R^{\rm ret}(z,\mathcal A_N)\|
    \le(\operatorname{Re} z)^{-1},
    \qquad \operatorname{Re} z>0.
\]
\end{proposition}

\begin{proposition}\label{prop:compact_case}
If \(\mathcal K\) is compact and \(\Pi_N\to I\) strongly, then
\[
    \|\Pi_N\mathcal K\Pi_N-\mathcal K\|\longrightarrow0.
\]
Consequently, Assumption~\ref{ass:stable_cvg} holds uniformly on compact subsets of \(\rho(\mathcal K)\). Assumption~\ref{ass:strong_stable_cvg} and Theorem~\ref{thm:hausdorff_cvg_isolated_eig} then apply to every isolated eigenvalue of finite algebraic multiplicity. In particular, they apply to every nonzero isolated eigenvalue of a compact operator.
\end{proposition}

\begin{remark}[Compact operators]
A stochastic Koopman or Markov operator with a square-integrable transition kernel is Hilbert--Schmidt and hence compact. Proposition~\ref{prop:compact_case} therefore applies when the stochastic transition kernel has this smoothing property.
\end{remark}

\begin{remark}\label{rm:scope_verifiable_conditions}
The spectral sets used in Section~\ref{sec:numerical_studies} may be compared with these ideal bounds. Proposition~\ref{prop:contractive_case} gives \(M_z=(R-1)^{-1}\) on \(\Gamma_R\) for a contractive discrete-time operator and \(R>1\). Proposition~\ref{prop:generator_measure_preserving} gives \(M_z=100\) on \(\operatorname{Re} z=\pm0.01\) for a measure-preserving generator. These estimates concern large-data Galerkin compressions in Gram-orthonormal coordinates. The explicit Euler pendulum map in Section~\ref{sec:pendulum} is not volume preserving. Its calculations on \(\Gamma_{1.01}\) therefore do not satisfy the stated contractive estimate. When stability has not been established, the ResDMD residual gives complementary information.
\end{remark}

\subsection{A Chatelin-type residual estimate for isolated eigenvalues}\label{chap:a_posteriori_bounds}
We record an operator-level residual estimate for eigenvalues obtained from finite-section methods such as EDMD \cite{Williams_2015}. The residual \(\eta_N\) contains the exact Koopman operator and is not generally available from snapshot data alone.

\begin{definition}[Arithmetic mean of approximate eigenvalues]\label{def:ari_mean_eval}
Let \(\lambda\) be an isolated eigenvalue of \(\mathcal K\) with algebraic multiplicity \(m_\lambda\), and let \(\Gamma_\lambda\) be an isolating contour as in Definition~\ref{def:spec_proj}. Suppose that Assumption~\ref{ass:stable_cvg} holds on \(\Gamma_\lambda\) and that Assumption~\ref{ass:strong_stable_cvg} holds relative to it. For all sufficiently large \(N\), let \(\{\mu_{N,j}\}_{j=1}^{m_\lambda}\) denote the eigenvalues of \(\mathcal K_N\) inside \(\Gamma_\lambda\), counted with algebraic multiplicity. Define
\begin{equation}
    \bar\lambda_N\coloneqq
    \frac{1}{m_\lambda}\sum_{j=1}^{m_\lambda}\mu_{N,j}.
\end{equation}
\end{definition}

\begin{definition}\label{def:comp_residual}
Let \(\mathcal P_N^\lambda\) be the Riesz projection associated with the finite-section eigenvalues inside \(\Gamma_\lambda\). Define
\[
    \eta_N\coloneqq
    \|(\mathcal K-\mathcal K_N)\mathcal P_N^\lambda\|.
\]
\end{definition}

\begin{proposition}[Residual estimate]\label{thm:a_posteriori_eigenvalue_bound}
Let \(\lambda\) be an isolated eigenvalue of \(\mathcal K\) with finite algebraic multiplicity, and let \(\Gamma_\lambda\) isolate \(\lambda\) from the remainder of \(\sigma(\mathcal K)\). Suppose that Assumption~\ref{ass:stable_cvg} holds on \(\Gamma_\lambda\) and that Assumption~\ref{ass:strong_stable_cvg} holds relative to it. Then, for all sufficiently large \(N\),
\begin{equation}\label{eq:chatelin_residual_bound}
    |\lambda-\bar\lambda_N|\le C\eta_N,
\end{equation}
where \(C\) is independent of \(N\).
\end{proposition}

The arithmetic mean is invariant under relabeling. It is also the quantity controlled by the trace when a multiple eigenvalue splits under finite-section approximation. If \(m_\lambda=1\), then \(\bar\lambda_N=\mu_N\), and Proposition~\ref{thm:a_posteriori_eigenvalue_bound} gives an estimate for the individual eigenvalue.

\begin{remark}
Section~\ref{sec:linear_oscillator} considers a finite-data analogue of Eq.~\eqref{eq:chatelin_residual_bound}. The regularized validation statistic used there is not an estimator of the operator quantity \(\eta_N\) without additional assumptions.
\end{remark}

%%%%%%%% Section 5 %%%%%%%%
\section{Numerical examples}\label{sec:numerical_studies}
The four examples address coefficient-coordinate sensitivity, finite-data residuals, contour geometry, and pseudomode-based signal separation. Unless stated otherwise, the reported quantities belong to the empirical finite section rather than to the limiting Koopman operator. Table~\ref{tab:numerical_settings} lists the data sets, dictionaries, and spectral evaluation sets.

\begin{table}[H]
\centering
\caption{Data sets, dictionaries, and spectral evaluation sets used in the numerical examples.}
\label{tab:numerical_settings}
\resizebox{\textwidth}{!}{%
\begin{tabular}{lll}
\toprule
Example & Data and dictionary & Spectral evaluation set \\
\midrule
Pendulum & 1500 snapshots, Fourier--Hermite dictionaries & circles and vertical lines \\
Linear oscillator & random Fourier features, assembly size \(30N\) & neighborhoods of known polynomial eigenvalues \\
Lorenz-type system & Hermite hyperbolic-cross dictionary, \(N=110\) & rectangular grid \\
Noisy oscillators & delay width \(L=850\), \(M_{\rm eff}=200\) & spectral grid for singular vectors \\
\bottomrule
\end{tabular}}
\end{table}

\subsection{Pendulum}\label{sec:pendulum}
The first example uses one explicit Euler step for an undamped pendulum with \(\Delta t=0.2\). Exact iterates of this discrete map provide the snapshot pairs, so no further time-integration error is introduced. The Euler map differs from the continuous Hamiltonian flow and is not volume preserving. With \(v=\dot\theta\) and \(c=3g/(2l)\), its Jacobian determinant is \(1+\Delta t^2c\cos\theta\). The continuous flow serves only as a qualitative reference. We retain raw singular values in nonorthogonal coefficient coordinates to examine their dependence on the threshold and dictionary size. This calculation motivates the Gram formulation in Section~\ref{sec:method}.

\begin{figure}
    \centering
    \includegraphics[width=0.95\linewidth]{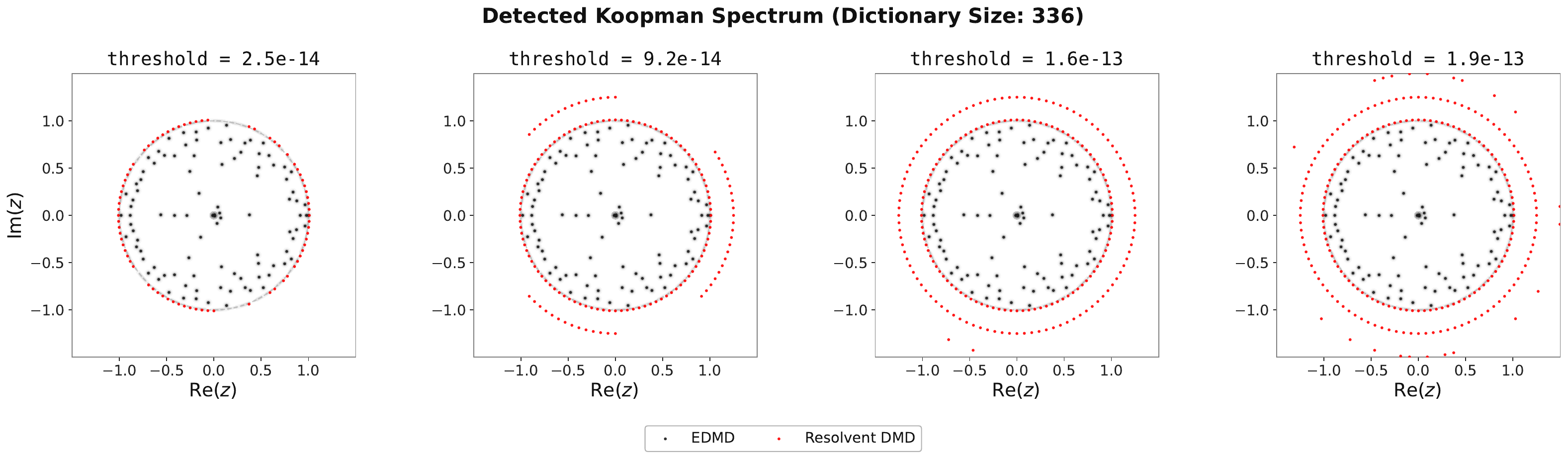}
    \caption{Coefficient-space threshold sensitivity for the Euler pendulum map with a 336-function subdictionary. Black points are EDMD eigenvalues. Red points are spectral samples satisfying the indicated inverse-resolvent threshold. All quantities refer to the finite matrix in the displayed coordinates.}
    \label{fig:pendulum_1}
\end{figure}

\begin{figure}
    \centering
    \includegraphics[width=0.95\linewidth]{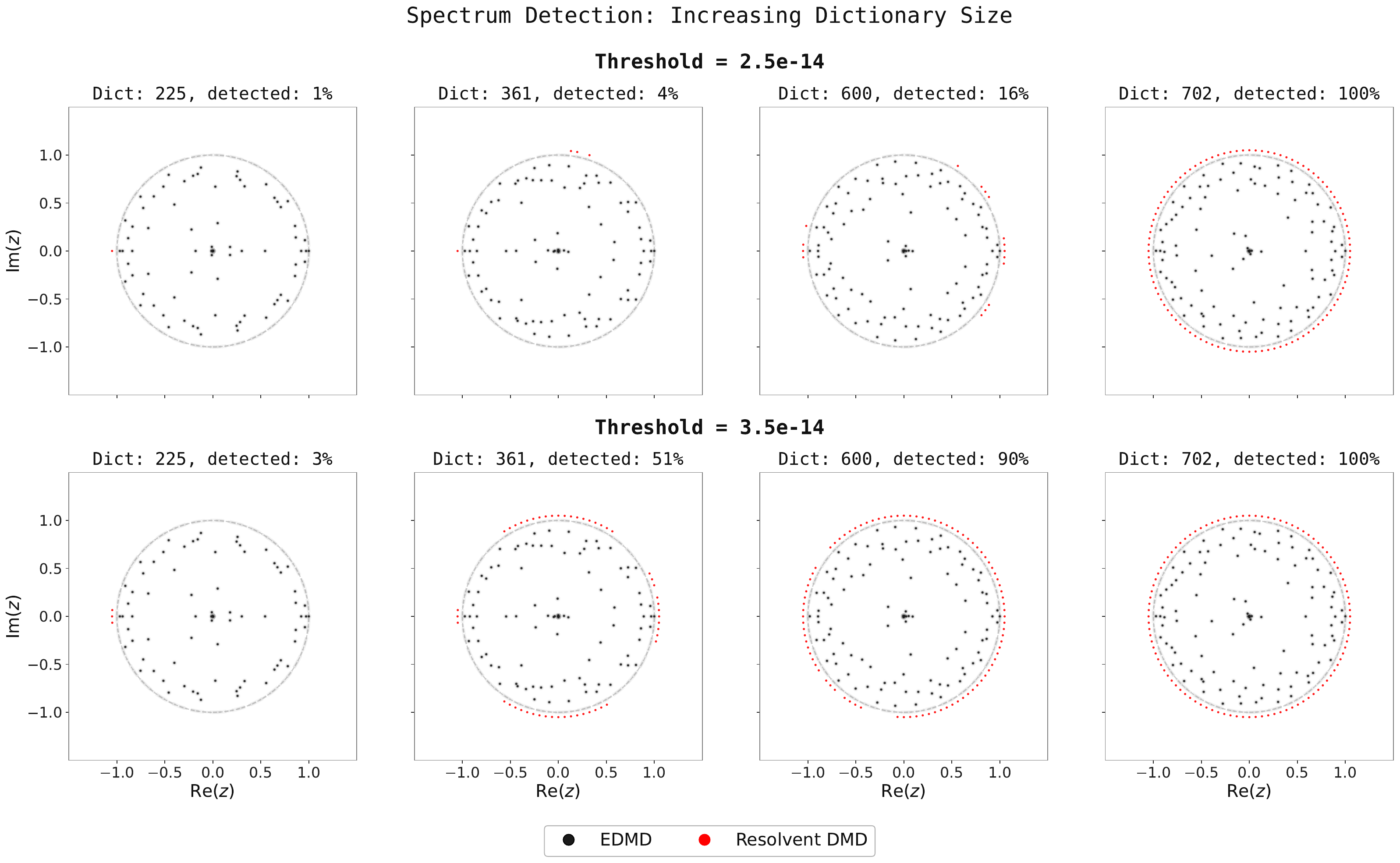}
    \caption{Fraction of spectral samples on \(\Gamma_{1.01}\) satisfying the coefficient-space inverse-resolvent threshold as the subdictionary size increases. The variation reflects changes in both the finite section and the coordinate conditioning.}
    \label{fig:pendulum_2}
\end{figure}

We generate one trajectory of 1500 steps. The initial angle is uniform on \([0,2\pi]\), and the initial angular velocity is uniform on \([-5,5]\). The continuous-time equation is
\[
    \ddot\theta=-\frac{3g}{2l}\sin\theta,
\]
with \(g=9.8\) and \(l=1.0\). The dictionary consists of products of Fourier modes \(\exp(in\theta)\), with \(n\in[-20,20]\), and Hermite polynomials in \(\dot\theta/\sqrt2\) through degree 30. The full candidate set contains \(41\times31=1271\) functions. Figure~\ref{fig:pendulum_1} uses 336 functions. Figure~\ref{fig:pendulum_2} uses nested subdictionaries of sizes \(225\), \(361\), \(600\), and \(702\), with the truncation and column order used in the implementation. These coordinates are not orthonormal for the empirical trajectory measure. The plotted Euclidean singular values therefore describe coefficient-space geometry rather than the Gram geometry of Proposition~\ref{prop:gram_identity}.

Figure~\ref{fig:pendulum_1} shows the spectral points that satisfy four coefficient-space inverse-resolvent thresholds. Figure~\ref{fig:pendulum_2} gives the selected fraction on \(\Gamma_{1.01}\) as the dictionary dimension increases. Both the finite matrix and the conditioning of its coordinate representation vary across these calculations. Remark~\ref{rm:coordinate_conditioning} bounds the discrepancy between coefficient and Gram geometries. The figures are therefore sensitivity studies and do not establish convergence of the Gram-weighted resolvent field. Appendix~\ref{app:pendulum_additional} contains additional generator and threshold calculations.

%%%%%%%% Section 4.2 %%%%%%%%
\subsection{Linear oscillator}\label{sec:linear_oscillator}
For the linear oscillator, polynomial Koopman eigenvalues are available in closed form. Consider
\[
    \dot x=A_{\rm sys}x,
    \qquad
    A_{\rm sys}=\begin{pmatrix}-0.3&-2\\2&-0.3\end{pmatrix},
\]
sampled at \(\Delta t=0.2\). The discrete flow map is \(T(x)=Ex\), where \(E=e^{A_{\rm sys}\Delta t}\). Polynomial observables have eigenvalues
\[
    \lambda_{n_1,n_2}=z_1^{n_1}z_2^{n_2},
    \qquad z_{1,2}=e^{(-0.3\pm2i)0.2},
    \qquad |z_1|\approx0.9418.
\]
These values serve as finite-data reference points. Whether they are isolated eigenvalues of an infinite-dimensional Koopman operator depends on the observable space and the sampling measure.

For each \(N\in\{20,30,50,80,120,180,250,350,500\}\), we assemble the EDMD matrix from \(M_{\rm train}=30N\) snapshot pairs. The dictionary consists of random Fourier features (RFFs) with kernel length scale \(\ell=1.5\) \cite{rahimi2007random}. The dictionaries are generated independently and are not nested. Thus the calculation compares a finite collection of approximation spaces rather than a Galerkin sequence. For a target \(\lambda_{n_1,n_2}\), let \(\mu_{N,M}\) be the nearest eigenvalue of \(\mathbf K_{N,M}\). Let \(v_{N,M}\) and \(w_{N,M}\) be the associated right and left eigenvectors, normalized by \(w_{N,M}^*v_{N,M}=1\). Define
\[
    g_{N,M}(x)=\sum_j(v_{N,M})_j\psi_j(x),
    \qquad
    r_{N,M}(x)=g_{N,M}(Ex)-\mu_{N,M}g_{N,M}(x).
\]

An independent sample of \(M_{\rm val}=10^5\) points from the same distribution is used to evaluate the residual norm and the Gram matrix \(\widehat G_{\rm val}\). We write \(Q_\alpha(\widehat G_{\rm val})\) for the chosen regularized inverse. It may be an inverse on a retained subspace, a truncated pseudoinverse, or a Tikhonov filter. The reported statistic is
\begin{equation}\label{eq:eta-simple}
  \widehat\eta_{N,M}^{\rm reg}
  =\left(w_{N,M}^* Q_\alpha(\widehat G_{\rm val})w_{N,M}\right)^{1/2}
  \left(\frac{1}{M_{\rm val}}\sum_{q=1}^{M_{\rm val}}|r_{N,M}(x_q^{\rm val})|^2\right)^{1/2}.
\end{equation}
The residual uses the exact flow map \(E\), so this evaluation introduces no temporal discretization error. Equation~\eqref{eq:eta-simple} combines an empirical Gram matrix with numerical regularization. It is therefore not the population residual \(\eta_N\) from Definition~\ref{def:comp_residual}. Five polynomial targets are matched independently for each \(N\).

\begin{figure}[ht]
\centering
\includegraphics[width=0.65\textwidth]{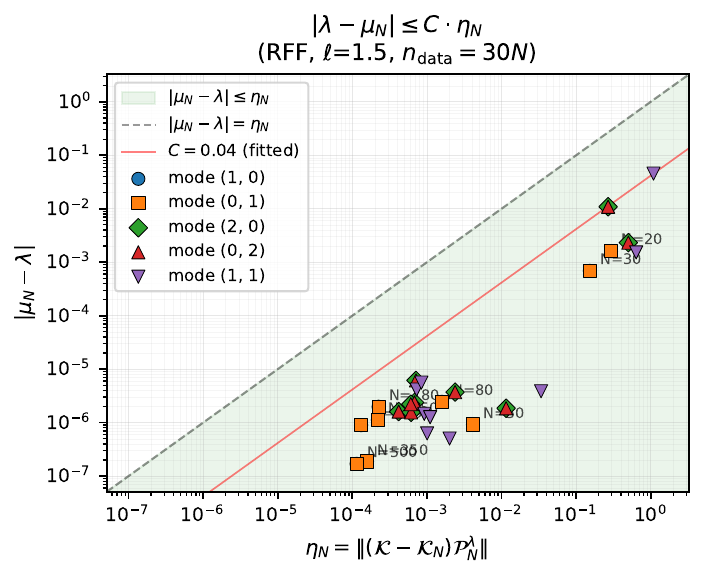}
\caption{Empirical eigenvalue error \(|\mu_{N,M}-\lambda|\) versus the regularized evaluation statistic \(\widehat\eta_{N,M}^{\rm reg}\) for five polynomial targets and nine RFF dictionary dimensions. The line \(C\approx0.04\) is fitted to these finite-data values and is not a theorem constant.}
\label{fig:thm54-scatter}
\end{figure}

\begin{figure}[ht]
\centering
\includegraphics[width=0.85\textwidth]{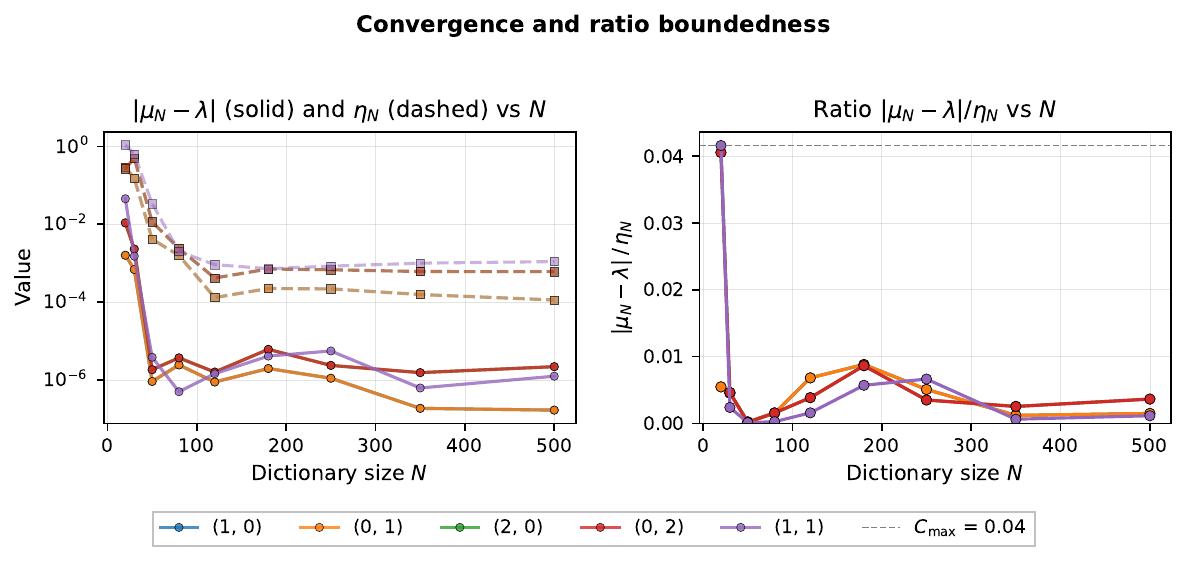}
\caption{Eigenvalue errors and regularized evaluation statistics across the tested RFF dimensions. The left panel shows both quantities. The right panel shows their ratio, and the dashed gray line marks the observed maximum \(0.042\).}
\label{fig:thm54-convergence}
\end{figure}

All 45 values in Figure~\ref{fig:thm54-scatter} lie below the diagonal. The largest observed ratio is
\[
  \max_{N,\lambda}\frac{|\mu_{N,M}-\lambda|}{\widehat\eta_{N,M}^{\rm reg}}=0.042.
\]
The line labeled \(C\approx0.04\) is fitted to the finite-data values and is not the constant in Proposition~\ref{thm:a_posteriori_eigenvalue_bound}. Figure~\ref{fig:thm54-convergence} shows the errors and residual statistics as functions of \(N\). The dictionaries are nonnested, and the required uniform stability assumptions are not verified. The plots therefore do not provide an asymptotic test of Proposition~\ref{thm:a_posteriori_eigenvalue_bound}.

%%%%%%%% Section 4.3 %%%%%%%%
\subsection{Lorenz-type system}\label{sec:comparing_pseudospectra}
\begin{figure}[!htb]
    \centering
    \includegraphics[width=\textwidth]{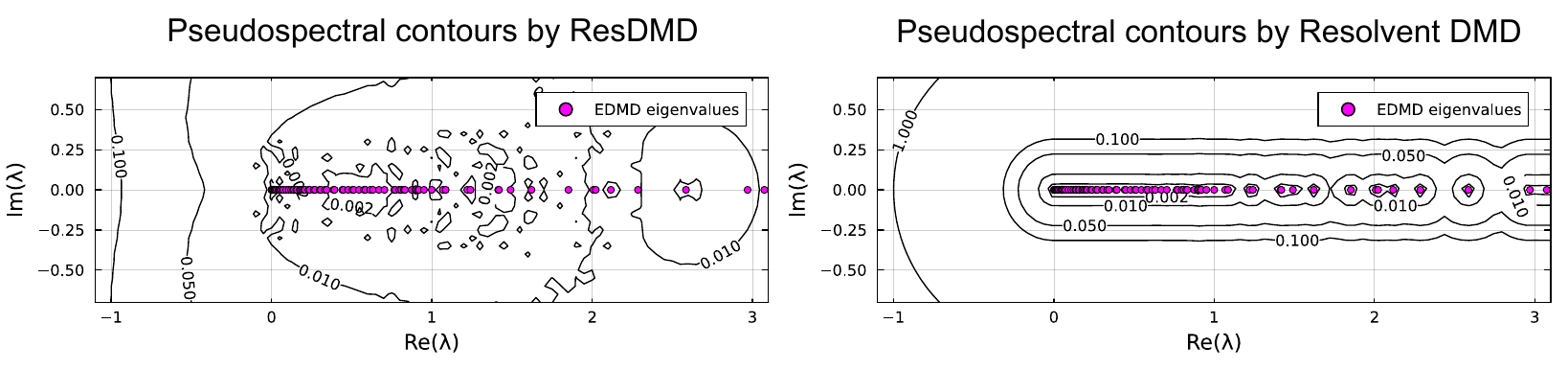}
    \caption{Contour fields for the modified Lorenz system. The left panel shows the ResDMD residual field. The right panel shows the reciprocal resolvent norm of the discrete-time EDMD finite section. Pink points denote EDMD eigenvalues. The comparison is qualitative because the plotted quantities are evaluated in coefficient coordinates.}
    \label{fig:pseudospectra_comparison}
\end{figure}

The Lorenz calculation compares a ResDMD residual field with the reciprocal resolvent norm of an EDMD finite section. Both fields use the same data set and dictionary. Proposition~\ref{prop:resdmd_gap} gives a sublevel-set inclusion when the quantities are unregularized and use the same retained Gram geometry. Figure~\ref{fig:pseudospectra_comparison} is instead computed in coefficient coordinates. We therefore compare only the contour geometry.

The modified Lorenz system is
\begin{equation}\label{eq:lorenz_system}
\begin{aligned}
    \dot x_1&=\sigma_{\rm L}(x_2-x_1),\\
    \dot x_2&=x_1(\rho-\kappa_{\rm L}x_3)-x_2,\\
    \dot x_3&=\kappa_{\rm L}x_1x_2-\beta x_3,
\end{aligned}
\end{equation}
with \(\sigma_{\rm L}=10\), \(\beta=8/3\), \(\rho=40\), and \(\kappa_{\rm L}=0.1\). We evaluate the time-\(\Delta t\) map numerically with \(\Delta t=0.05\). A hyperbolic-cross truncation of tensor products of normalized Hermite polynomials gives \(N=110\). The same Gauss--Hermite rule is used for the coefficient inner products and the matrix assembly. Its quadrature measure is not assumed to be invariant under the Lorenz flow.

Figure~\ref{fig:pseudospectra_comparison} shows the ResDMD residual contours and the reciprocal finite-section resolvent norm. The resolvent contours are smoother on the chosen grid. The calculations use coefficient coordinates and different discrete formulas, so the visual difference does not establish a difference in accuracy.

%%%%%%%% Section 4.4 %%%%%%%%
\subsection{Noisy two-oscillator signal}\label{sec:noisy_oscillators}
The final example considers a noisy scalar signal generated by two weakly damped oscillators. The latent four-dimensional process is Markovian, but a finite delay vector from one scalar observation need not be a Markov state. We therefore regard the matrix in Eq.~\eqref{eq:oscillator_galerkin} as a regularized one-step evolution matrix on the selected delay space. It is not treated as an exact finite-dimensional Koopman representation.

Let \((X_i,V_i)\) denote the position and velocity of oscillator \(i\).
\begin{align}
    dX_i(t)&=V_i(t)\,dt,\notag\\
    dV_i(t)&=\bigl(-2\zeta\omega_iV_i(t)-\omega_i^2X_i(t)\bigr)\,dt
              +\sigma_i\,dW_i(t),\qquad i=1,2,
\end{align}
where \(\omega_1=3\pi/5\), \(\omega_2=6\pi/5=2\omega_1\), \(\zeta=0.2\), \(\sigma_1=0.05\), and \(\sigma_2=0.15\). The observed scalar signal is
\[
    Y(t)=X_1(t)+X_2(t)+\xi(t),
\]
where the sampled measurement errors are independent \(N(0,0.015^2)\) variables and are independent of the Wiener processes \(W_1\) and \(W_2\).

We use 1050 observations with sampling interval \(\Delta t_{\rm obs}=0.1\). Delay-coordinate constructions of this form are standard in Koopman and DMD analysis \cite{arbabi2017ergodic}. A delay width of \(L=850\) gives \(M_{\rm eff}=200\) consecutive delay-vector pairs. With uniform weights,
\begin{equation}\label{eq:oscillator_galerkin}
    \mathbf K_{L,M_{\rm eff}}
      =(\Psi_X^*W\Psi_X)^\dagger(\Psi_X^*W\Psi_Y).
\end{equation}
Since \(\operatorname{rank}(\Psi_X^*W\Psi_X)\le200<850\), the SVD cutoff determines the retained space for every shifted problem. This is the regime with fewer snapshots than dictionary elements considered for ResDMD in \cite{colbrook2024fewer}.

Let \(X\) and \(Y\) denote the current and successor delay vectors. The conditional decomposition
\begin{equation}\label{eq:conditional_variance}
 \mathbb E\!\left[|g(Y)-zg(X)|^2\mid X\right]
 =|\mathbb E[g(Y)\mid X]-zg(X)|^2
  +\operatorname{Var}(g(Y)\mid X)
\end{equation}
shows that one observed successor contains conditional variance in addition to the squared Markov-operator residual. A deterministic ResDMD residual therefore requires variance information or replicated successors \cite{colbrook2024beyond}. We use a residual-based ResDMD surrogate and compare it with the pseudomodes of the retained delay-space evolution matrix.

\begin{figure}[!htb]
    \centering
    \includegraphics[width=0.49\textwidth]{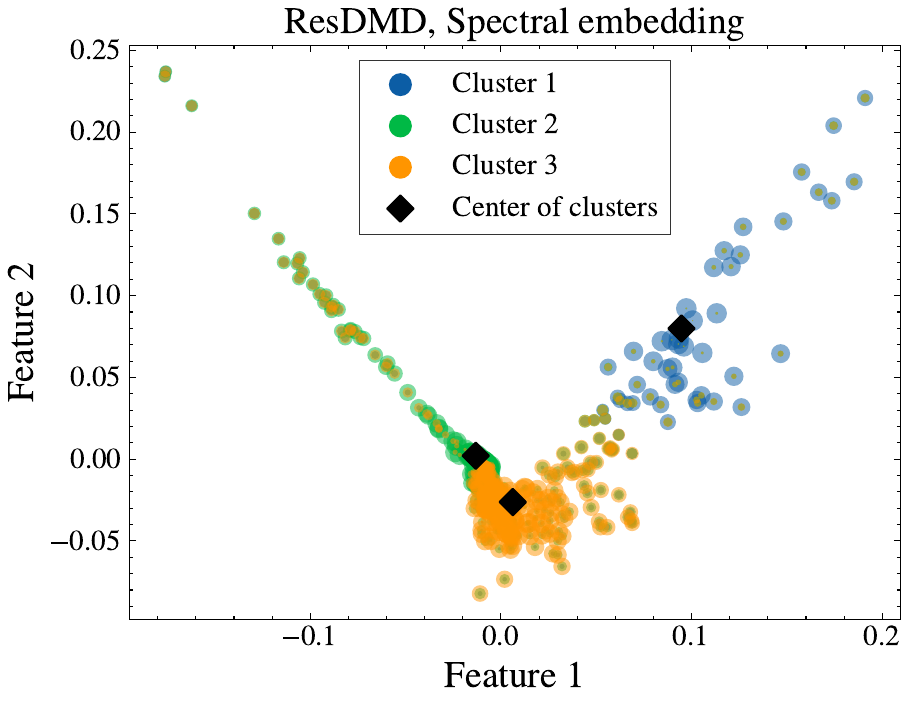}\hfill
    \includegraphics[width=0.49\textwidth]{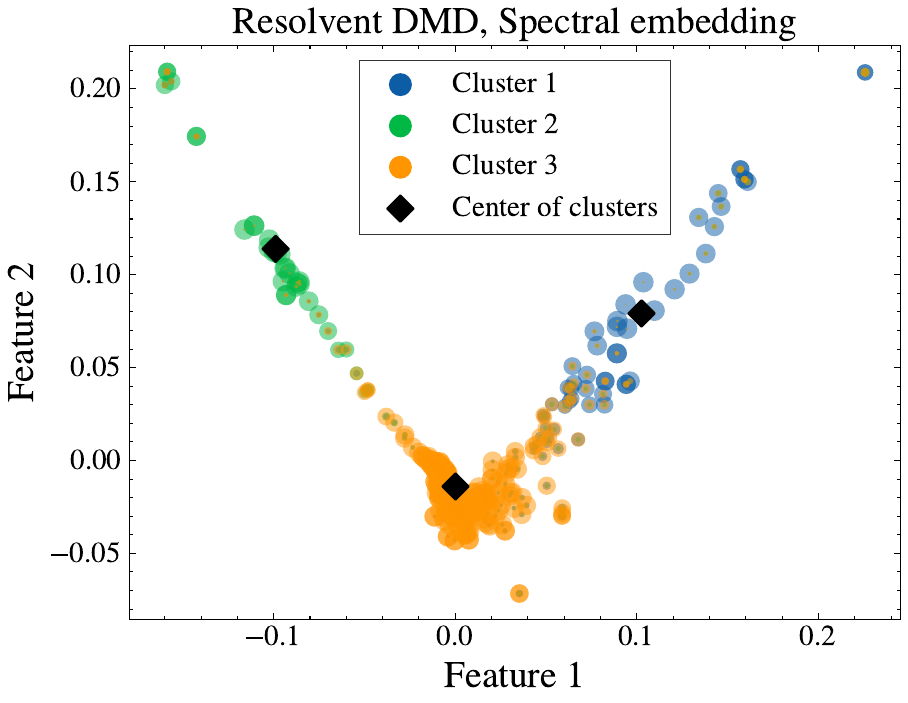}
    \caption{Spectral embeddings constructed from the ResDMD residual surrogate on the left and the Resolvent DMD pseudomodes on the right. The three groups overlap less in the right panel for this realization.}
    \label{fig:embedding_two_osc}
\end{figure}

\begin{figure}[!htb]
    \centering
    \includegraphics[width=0.49\textwidth]{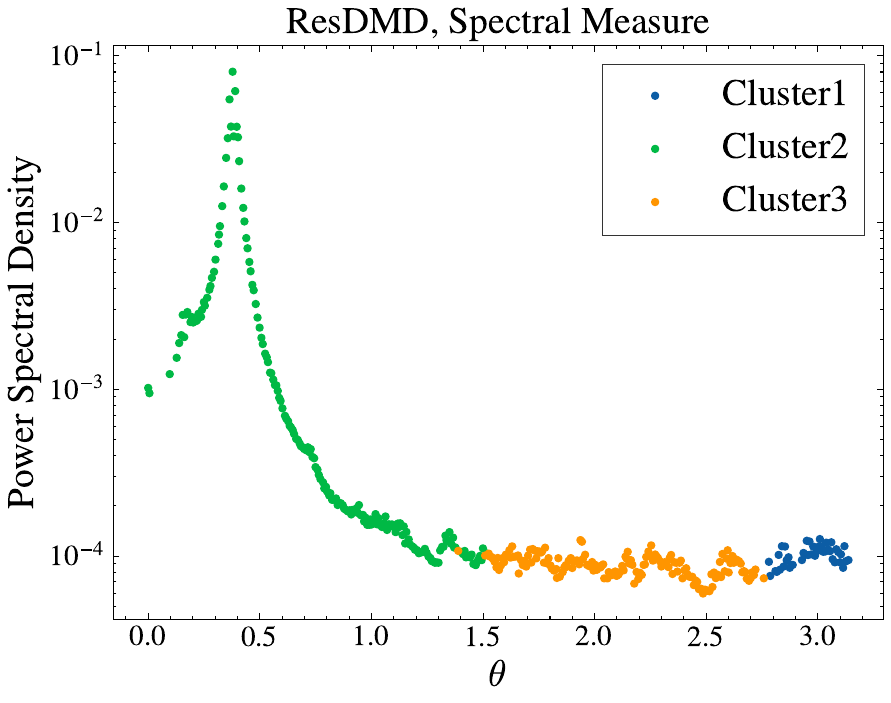}\hfill
    \includegraphics[width=0.49\textwidth]{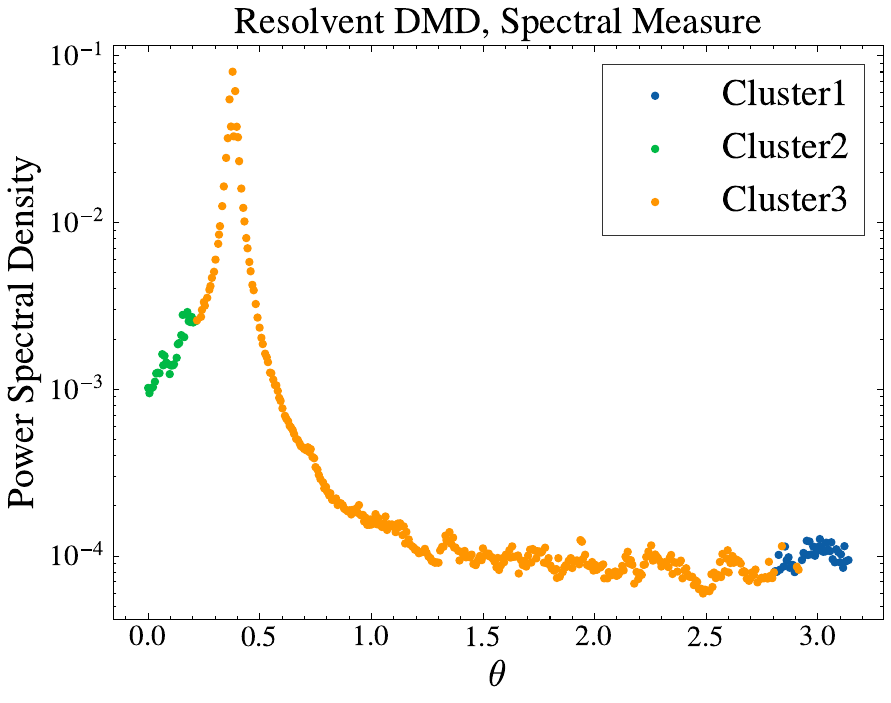}
    \caption{Cluster-labeled empirical power spectral densities for the ResDMD observables on the left and the Resolvent DMD pseudomodes on the right. The two principal frequency bands are more clearly separated in the right panel for this realization.}
    \label{fig:spectral_two_osc}
\end{figure}

\begin{figure}[!htb]
    \centering
    \includegraphics[width=0.49\textwidth]{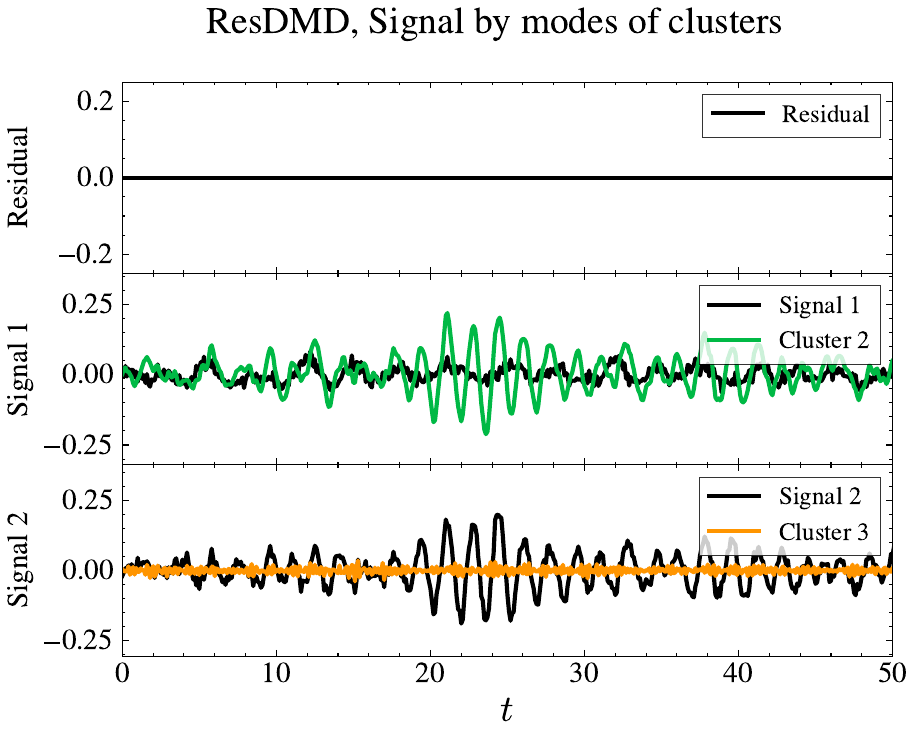}\hfill
    \includegraphics[width=0.49\textwidth]{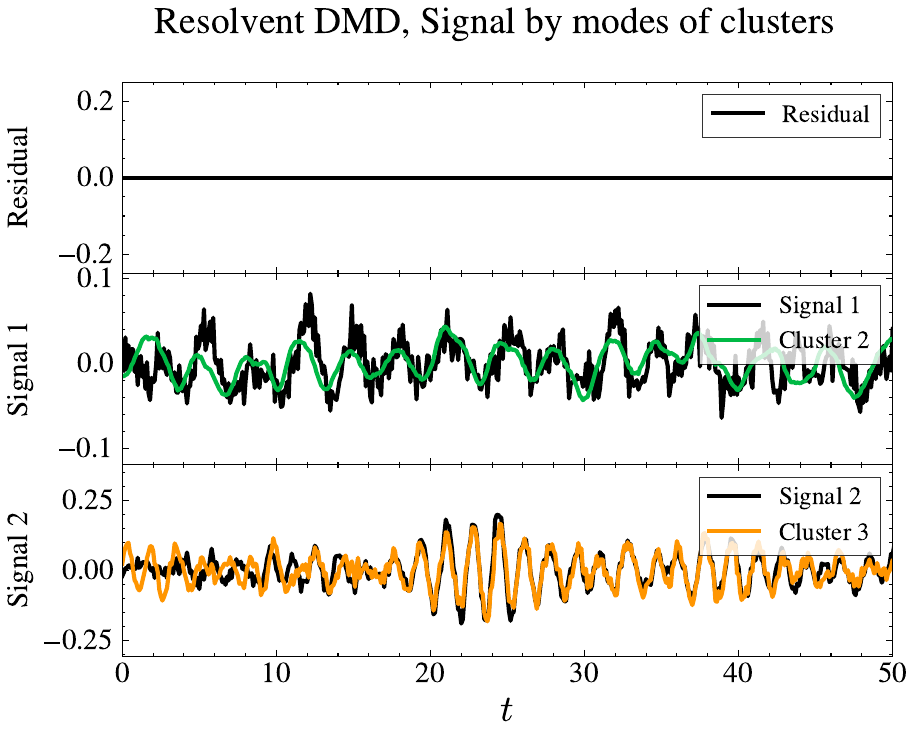}
    \caption{Time-domain contributions reconstructed from the clustered observables. The left panel uses the ResDMD surrogate, and the right panel uses the Resolvent DMD pseudomodes. In this realization, the right-hand reconstruction follows the two simulated components more closely.}
    \label{fig:reconst_two_osc}
\end{figure}

At each spectral point \(z\), both calculations return an observable. Following \cite{sakata2024enhancing}, we form principal-angle similarities and then apply spectral embedding and fuzzy \(c\)-means clustering. The resulting groups align with the two simulated frequency bands and a noise-dominated component. In the realization considered here, the Resolvent DMD embedding has less cluster overlap than the ResDMD surrogate in Figure~\ref{fig:embedding_two_osc}. The cluster-labeled power spectral densities in Figure~\ref{fig:spectral_two_osc} are empirical spectra of the computed observables. They are not Koopman spectral measures in the sense of \cite{korda2020data}. Figure~\ref{fig:reconst_two_osc} shows the corresponding time-domain contributions. These comparisons are qualitative and concern one realization.

%%%%%%%% Section 6 %%%%%%%%
\section{Discussion}\label{sec:discussion}
The finite-dimensional object considered here is the Koopman representation in the empirical observable geometry. After rank truncation and Gram orthonormalization, the inverse of the smallest shifted singular value is the resolvent norm of the retained matrix. The corresponding singular vectors give a minimum-residual pseudomode and its optimal forcing direction. Their computation does not require formation of the shifted inverse.

Proposition~\ref{prop:resdmd_gap} separates the deterministic ResDMD residual into two terms. The shifted finite-section residual measures the projected pseudoeigenvalue error. The invariance defect measures the component of the propagated observable outside the retained trial space. Both terms are needed for the full residual. The distinction is important when the dictionary is not invariant or the empirical Gram matrix is rank truncated.

Convergence of the Galerkin projections alone does not give an operator-level interpretation. Pointwise stability yields retained strong resolvent convergence and one-sided inclusion of inverse-resolvent sublevel sets. Stability on a punctured isolating neighborhood gives local Hausdorff convergence near an isolated eigenvalue. The correct algebraic multiplicity additionally requires stability of the Riesz projection ranks. Contractivity and skew-adjointness provide explicit resolvent bounds, while compactness gives norm convergence. The backward shift shows why finite-section eigenvalues alone may miss spectral information carried by the resolvent.

The numerical examples have distinct and limited purposes. The pendulum and Lorenz calculations use coefficient coordinates, so they measure sensitivity rather than Gram-weighted convergence. The linear-oscillator statistic is a finite-data analogue of the operator residual. The noisy-oscillator comparison uses one stochastic realization. Further analysis is needed for empirical estimation of the invariance defect, variance-corrected stochastic residuals, and selection of dictionaries, retained ranks, and spectral evaluation sets.

%%%%%%%% Section 7 %%%%%%%%
\section{Conclusion}\label{sec:conclusion}
The smallest singular triplet of a Gram-orthonormalized Koopman finite section describes resolvent amplification in the retained observable norm. Its singular value is the reciprocal resolvent norm, and its singular vectors determine a pseudomode and the associated forcing direction. On the same retained space, the deterministic ResDMD residual splits into a shifted finite-section residual and an invariance defect. Finite-section stability then gives one-sided inclusion of inverse-resolvent sublevel sets and local spectral convergence near isolated eigenvalues. These conclusions first apply to the retained finite section. Their extension to the Koopman operator requires the stability assumptions stated above.

\appendix
\section{Technical details and proofs}\label{app:proofs}

\subsection{Finite-data and retained-resolvent preliminaries}\label{app:finite_data_limit}
The large-data limit is taken with \(N\) fixed and \(M\to\infty\). Under standard independent-sampling or ergodic assumptions, the empirical Gram and cross-Gram matrices converge almost surely to their population limits \cite{korda2018convergence,Williams_2015}. A positive-definite population Gram matrix implies eventual nonsingularity of the empirical matrices. For a rank-deficient dictionary, one may instead fix the retained rank and assume a spectral gap at the truncation threshold. The corresponding truncated pseudoinverse is then used. Under either condition,
\[
  \mathbf K_{N,M}\longrightarrow\mathbf K_N,
  \qquad
  (zI-\mathbf K_{N,M})^{-1}\longrightarrow(zI-\mathbf K_N)^{-1}
\]
locally uniformly on compact subsets of \(\rho(\mathbf K_N)\). Convergence of the Gram blocks alone is not sufficient when the retained rank is unstable. In that case, the pseudoinverse need not converge.

For the ideal Galerkin compressions, \(\Pi_N\mathcal K\Pi_N\to\mathcal K\) strongly. At a point satisfying Assumption~\ref{ass:stable_cvg},
\[
\mathcal R_N^{\rm ret}(z)\Pi_N f-\Pi_N\mathcal R(z)f
=-\mathcal R_N^{\rm ret}(z)\Pi_N\mathcal K(I-\Pi_N)\mathcal R(z)f.
\]
The right-hand side converges to zero by the uniform resolvent bound and the strong convergence of \(\Pi_N\). Since \(\Pi_N\mathcal R(z)f\to\mathcal R(z)f\), Eq.~\eqref{eq:resolvent_strong_cvg} follows.

\subsection*{Proof of Proposition~\ref{prop:gram_identity}}
\begin{proof}
Since the columns of \(V_r\) are orthonormal eigenvectors associated with the retained spectrum of \(G_{N,M}\) and \(V_r^*G_{\rm disc}V_r=0\),
\[
  \Phi_{X,r}^*W\Phi_{X,r}
  =\Sigma_r^{-1/2}V_r^*G_{N,M}V_r\Sigma_r^{-1/2}
  =I_r.
\]
Hence the Euclidean norm in retained coordinates agrees with the empirical observable norm. Moreover,
\[
\begin{aligned}
  \Phi_{X,r}^*W\Psi_YV_r\Sigma_r^{-1/2}
  &=\Sigma_r^{-1/2}V_r^*B_{N,M}V_r\Sigma_r^{-1/2}\\
  &=\widetilde K_{N,M},
\end{aligned}
\]
so \(\widetilde K_{N,M}\) is the least-squares compression in this basis. Equation~\eqref{eq:gram_identity} is the standard inverse-norm identity for the smallest singular value. If \(G_{N,M}=V\Sigma V^*\) is nonsingular, then
\[
G_{N,M}^{1/2}\mathbf K_{N,M}G_{N,M}^{-1/2}
=V\widetilde K_{N,M}V^*,
\]
and the corresponding shifted matrices have the same singular values. Finally, any two orthonormal bases of the retained empirical subspace differ by a unitary transformation.
\end{proof}

\subsection*{Proof of Proposition~\ref{prop:resdmd_gap}}
\begin{proof}
Set \(\Phi=\Phi_{X,r}\). For \(c=V_r\Sigma_r^{-1/2}q\),
\[
  \Psi_Xc=\Phi q,
  \qquad
  P_r\Psi_Yc
  =\Phi\Phi^*W\Psi_YV_r\Sigma_r^{-1/2}q
  =\Phi\widetilde K_{N,M}q.
\]
Therefore,
\[
  \Psi_Yc-z\Psi_Xc
  =\Phi(\widetilde K_{N,M}-zI_r)q
   +(I-P_r)\Psi_Yc.
\]
The two terms are \(W\)-orthogonal. Since \(\Phi^*W\Phi=I_r\), the Pythagorean identity gives Eq.~\eqref{eq:resdmd_decomposition}. The same identity shows that \(\|q\|_2=1\) if and only if \(\|\Psi_Xc\|_W=1\). Minimization over unit vectors gives \(\tau_{N,M}(z)\ge s_{N,M}(z)\). The sublevel-set inclusion follows.
\end{proof}

\subsection*{Proof of Lemma~\ref{lem:uniform_bound}}
\begin{proof}
Fix \(z\in U\). By Assumption~\ref{ass:stable_cvg}, there exist \(N(z)\) and \(M_z<\infty\) such that, for every \(N>N(z)\),
\[
    z\in\rho(\mathcal K_N|_{\mathcal F_N}),
    \qquad
    \|\mathcal R_N(z)\|\le M_z.
\]
Set \(r_z=(2M_z)^{-1}\).

Let \(|w-z|<r_z\) and \(N>N(z)\). Since invertibility at \(w\) is not yet known, introduce the Neumann series
\[
S_N(w)\coloneqq\sum_{k=0}^{\infty}(z-w)^k\mathcal R_N(z)^{k+1}.
\]
It converges in operator norm because \(|w-z|\,\|\mathcal R_N(z)\|<1/2\), and
\[
S_N(w)=\bigl[I+(w-z)\mathcal R_N(z)\bigr]^{-1}\mathcal R_N(z).
\]
The factorization
\[
    wI-\mathcal K_N
    =(zI-\mathcal K_N)\bigl[I+(w-z)\mathcal R_N(z)\bigr]
\]
shows that \(S_N(w)=(wI-\mathcal K_N)^{-1}=\mathcal R_N(w)\). Moreover,
\[
\|\mathcal R_N(w)\|
\le
\frac{\|\mathcal R_N(z)\|}{1-|w-z|\,\|\mathcal R_N(z)\|}
\le 2M_z.
\]

The balls \(\{B(z,r_z)\mid z\in U\}\) form an open cover of the compact set \(U\). Choose a finite subcover \(\{B(z_j,r_{z_j})\}_{j=1}^J\), and set
\[
N_U\coloneqq\max_{1\le j\le J}N(z_j),
\qquad
M_U\coloneqq2\max_{1\le j\le J}M_{z_j}.
\]
For \(N>N_U\) and \(w\in U\), choose \(j\) such that \(w\in B(z_j,r_{z_j})\). Then \(w\in\rho(\mathcal K_N|_{\mathcal F_N})\) and \(\|\mathcal R_N(w)\|\le M_U\).
\end{proof}

\subsection*{Proof of Corollary~\ref{cor:spectral_projection}}
\begin{proof}
For \(f\in\mathcal F\), Definition~\ref{def:spec_proj} gives
\[
(\mathcal P_N^{\lambda}-\mathcal P^{\lambda})f
=
\frac{1}{2\pi i}\int_{\Gamma_\lambda}
\bigl(\mathcal R_N^{\rm ret}(z)\Pi_N f-\mathcal R(z)f\bigr)\,dz.
\]

For each \(z\in\Gamma_\lambda\), Eq.~\eqref{eq:resolvent_strong_cvg} gives \(\mathcal R_N^{\rm ret}(z)\Pi_N f\to\mathcal R(z)f\). Lemma~\ref{lem:uniform_bound}, applied to \(\Gamma_\lambda\), yields constants \(N_\Gamma\) and \(M_\Gamma\) such that
\[
\sup_{z\in\Gamma_\lambda,\,N>N_\Gamma}
\|\mathcal R_N^{\rm ret}(z)\|\le M_\Gamma.
\]
Since \(\mathcal R(\cdot)\) is continuous on the compact contour, the integrand is dominated by
\[
\bigl\|\mathcal R_N^{\rm ret}(z)\Pi_N f-\mathcal R(z)f\bigr\|
\le
\left(M_\Gamma+
\sup_{w\in\Gamma_\lambda}\|\mathcal R(w)\|\right)\|f\|.
\]
Dominated convergence for contour integrals therefore gives
\[
\|(\mathcal P_N^{\lambda}-\mathcal P^{\lambda})f\|
\le
\frac{1}{2\pi}\int_{\Gamma_\lambda}
\bigl\|\mathcal R_N^{\rm ret}(z)\Pi_N f-\mathcal R(z)f\bigr\|\,|dz|
\longrightarrow0.
\]
Thus \(\mathcal P_N^\lambda\to\mathcal P^\lambda\) strongly. By \cite[Proposition~3.10]{chaitin1983spectral}, strong convergence of these projections implies
\[
    \dim(\mathcal P_N^\lambda\mathcal F)
    \ge\dim(\mathcal P^\lambda\mathcal F)
\]
for all sufficiently large \(N\).
\end{proof}

\subsection*{Proof of Theorem~\ref{thm:hausdorff_cvg_isolated_eig}}
\begin{proof}
We first prove the two directed Hausdorff bounds and then address multiplicity.

\textit{Lower bound.}
Choose \(\varepsilon>0\) such that
\[
\gamma_\varepsilon\coloneqq\{z\mid |z-\lambda|=\varepsilon\}
\]
lies in \(\Delta\cap\rho(\mathcal K)\) and encloses no spectral point other than \(\lambda\). By hypothesis, the finite sections are stable on \(\gamma_\varepsilon\). Let \(\mathcal P^{\gamma_\varepsilon}\) and \(\mathcal P_N^{\gamma_\varepsilon}\) denote the corresponding Riesz projections. Since
\[
    \dim(\mathcal P^{\gamma_\varepsilon}\mathcal F)=m,
\]
Corollary~\ref{cor:spectral_projection} gives
\[
    \dim(\mathcal P_N^{\gamma_\varepsilon}\mathcal F)\ge m
\]
for all sufficiently large \(N\). Hence \(\mathcal P_N^{\gamma_\varepsilon}\ne0\), so
\[
    \sigma(\mathcal K_N)\cap B(\lambda,\varepsilon)\ne\varnothing.
\]
Equivalently,
\[
    \operatorname{dist}\bigl(\lambda,\sigma(\mathcal K_N)\cap\Delta\bigr)
    <\varepsilon
\]
for all sufficiently large \(N\).

\textit{Upper bound.}
Define the compact set
\[
    \overline\Delta_\varepsilon
    \coloneqq\overline\Delta\setminus B(\lambda,\varepsilon).
\]
Because \(\sigma(\mathcal K)\cap\overline\Delta_\varepsilon=\varnothing\), we have \(\overline\Delta_\varepsilon\subset\rho(\mathcal K)\). Lemma~\ref{lem:uniform_bound} therefore implies
\[
    \overline\Delta_\varepsilon
    \subset\rho(\mathcal K_N|_{\mathcal F_N})
\]
for all sufficiently large \(N\). Hence every \(\mu\in\sigma(\mathcal K_N)\cap\Delta\) satisfies \(|\mu-\lambda|<\varepsilon\), and
\[
    \sup_{\mu\in\sigma(\mathcal K_N)\cap\Delta}
    \operatorname{dist}(\mu,\{\lambda\})<\varepsilon.
\]
Combining the two bounds and letting \(\varepsilon\downarrow0\) proves the Hausdorff convergence.

\textit{Multiplicity.}
Suppose that Assumption~\ref{ass:strong_stable_cvg} holds relative to \(\Gamma\). Then \(\dim(\mathcal P_N^\lambda\mathcal F)=m\) for all sufficiently large \(N\). By the Riesz decomposition, the spectrum of \(\mathcal K_N\) on \(\mathcal P_N^\lambda\mathcal F\) is \(\sigma(\mathcal K_N)\cap\Delta\) \cite[Theorem 2.27]{chaitin1983spectral}. This invariant subspace has dimension \(m\). Hence the eigenvalues in \(\Delta\) have total algebraic multiplicity \(m\).
\end{proof}

\subsection*{Proof of Theorem~\ref{thm:one_sided_indicator}}
\begin{proof}
For the pointwise statement, let \(\delta>0\) and choose a unit vector \(f\in\mathcal F\) such that
\[
    \|\mathcal R(z)f\|\ge\|\mathcal R(z)\|-\delta.
\]
Equation~\eqref{eq:resolvent_strong_cvg} gives
\[
    \mathcal R_N^{\rm ret}(z)\Pi_N f
    \longrightarrow\mathcal R(z)f.
\]
Since
\[
    \|\mathcal R_N^{\rm ret}(z)\Pi_N f\|
    \le\|\mathcal R_N(z)\|,
\]
it follows that
\[
    \liminf_{N\to\infty}\|\mathcal R_N(z)\|
    \ge\|\mathcal R(z)f\|
    \ge\|\mathcal R(z)\|-\delta.
\]
Letting \(\delta\downarrow0\) proves the pointwise inequality and hence the pointwise inverse-resolvent inclusion.

For the uniform statement, set
\[
    a\coloneqq\inf_{z\in Q}\|\mathcal R(z)\|>\varepsilon^{-1}
\]
and choose \(\eta>0\) such that \(a-3\eta>\varepsilon^{-1}\). For each \(z_0\in Q\), choose a unit vector \(f_{z_0}\) with
\[
    \|\mathcal R(z_0)f_{z_0}\|>a-\eta.
\]
Continuity of \(w\mapsto\mathcal R(w)f_{z_0}\) gives \(r_{z_0}>0\) such that
\[
    \|\mathcal R(w)f_{z_0}\|>a-2\eta
\]
whenever \(|w-z_0|\le r_{z_0}\).

Let \(Q_{z_0}=Q\cap\overline{B(z_0,r_{z_0})}\). Lemma~\ref{lem:uniform_bound} bounds \(\mathcal R_N^{\rm ret}(w)\) uniformly on \(Q_{z_0}\), and the set
\[
    \{\mathcal R(w)f_{z_0}\mid w\in Q_{z_0}\}
\]
is compact. Since strong convergence of \(\Pi_N\) is uniform on compact subsets of \(\mathcal F\), the retained resolvent identity yields
\[
\begin{aligned}
&\sup_{w\in Q_{z_0}}
\|\mathcal R_N^{\rm ret}(w)\Pi_N f_{z_0}
  -\mathcal R(w)f_{z_0}\| \\
&\quad\le
\sup_{w\in Q_{z_0}}\|\mathcal R_N^{\rm ret}(w)\|\,\|\mathcal K\|
\sup_{w\in Q_{z_0}}\|(I-\Pi_N)\mathcal R(w)f_{z_0}\|\\
&\qquad+
\sup_{w\in Q_{z_0}}\|(I-\Pi_N)\mathcal R(w)f_{z_0}\|
\longrightarrow0.
\end{aligned}
\]
Therefore \(\|\mathcal R_N(w)\|>a-3\eta>\varepsilon^{-1}\) on \(Q_{z_0}\) for all sufficiently large \(N\). A finite subcover of \(Q\) completes the proof.
\end{proof}

\subsection*{Proof of Proposition~\ref{prop:contractive_case}}
\begin{proof}
Since \(\Pi_N\) is an orthogonal projection, \(\|\mathcal K_N\|\le\|\mathcal K\|\le1\). For \(|z|>1\), the Neumann series converges and gives
\[
\begin{aligned}
    \|(zI-\mathcal K_N)^{-1}\|
    &\le |z|^{-1}\sum_{j=0}^{\infty}|z|^{-j}\|\mathcal K_N^j\|\\
    &\le |z|^{-1}\sum_{j=0}^{\infty}|z|^{-j}
    =\frac{1}{|z|-1}.
\end{aligned}
\]
\end{proof}

\subsection*{Proof of Proposition~\ref{prop:generator_measure_preserving}}
\begin{proof}
For \(f,g\in\mathcal F_N\), skew-adjointness of \(\mathcal A\) gives
\[
    \langle\mathcal A_N f,g\rangle
    =\langle\mathcal A f,g\rangle
    =-\langle f,\mathcal A g\rangle
    =-\langle f,\mathcal A_N g\rangle.
\]
Thus \(\mathcal A_N\) is skew-adjoint on \(\mathcal F_N\). It is therefore normal with spectrum on \(i\mathbb R\), and the resolvent formula in part (i) follows.

To prove part (ii), first take \(f\in\bigcup_N\mathcal F_N\) and set \(h=(zI-\mathcal A)f\). For all sufficiently large \(N\), \(f\in\mathcal F_N\) and
\[
    \Pi_N h
    =\Pi_N(zI-\mathcal A)f
    =(zI_{\mathcal F_N}-\mathcal A_N)f.
\]
Consequently,
\[
    \mathcal R_N^{\mathcal A}(z)h
    =\mathcal R^{\rm ret}(z,\mathcal A_N)\Pi_N h
    =f.
\]
Because \(\bigcup_N\mathcal F_N\) is a core for \(\mathcal A\), the set \((zI-\mathcal A)(\bigcup_N\mathcal F_N)\) is dense in \(L^2(\mu)\). The uniform bound in part (i) extends the convergence from this dense set to all \(h\in L^2(\mu)\).

For part (iii), let \(\lambda\in\sigma(\mathcal A)\) and set \(z=\lambda+\delta\) with \(\delta>0\). Part (ii) and lower semicontinuity of the operator norm under strong convergence imply
\[
    \liminf_{N\to\infty}\|\mathcal R_N^{\mathcal A}(z)\|
    \ge\|\mathcal R(z,\mathcal A)\|
    =\frac{1}{\operatorname{dist}(z,\sigma(\mathcal A))}
    \ge\frac{1}{\delta}.
\]
Since \(\mathcal A_N\) is normal,
\[
    \|\mathcal R_N^{\mathcal A}(z)\|
    =\frac{1}{\operatorname{dist}(z,\sigma(\mathcal A_N))}.
\]
Thus
\[
    \limsup_{N\to\infty}\operatorname{dist}(z,\sigma(\mathcal A_N))
    \le\delta,
\]
and the triangle inequality gives
\[
    \limsup_{N\to\infty}\operatorname{dist}(\lambda,\sigma(\mathcal A_N))
    \le2\delta.
\]
Letting \(\delta\downarrow0\) proves (iii). For a contraction semigroup, each compression \(\mathcal A_N\) is dissipative, and the stated estimate follows from the standard numerical-range resolvent bound for finite-dimensional dissipative operators \cite{engel1999one,pazy2012semigroups}.
\end{proof}

\subsection*{Proof of Proposition~\ref{prop:compact_case}}
\begin{proof}
Since \(\Pi_N\) is orthogonal,
\[
\begin{aligned}
    \|\mathcal K-\Pi_N\mathcal K\Pi_N\|
    &\le\|(I-\Pi_N)\mathcal K\|
       +\|\Pi_N\mathcal K(I-\Pi_N)\|\\
    &\le\|(I-\Pi_N)\mathcal K\|
       +\|\mathcal K(I-\Pi_N)\|.
\end{aligned}
\]
The first term converges to zero because \(\mathcal K\) is compact and \(\Pi_N\to I\) strongly. The second term equals \(\|(I-\Pi_N)\mathcal K^*\|\) and also converges to zero, since \(\mathcal K^*\) is compact. Hence \(\mathcal K_N\to\mathcal K\) in operator norm. The second resolvent identity gives uniform resolvent convergence on compact subsets of \(\rho(\mathcal K)\). For an isolated eigenvalue of finite algebraic multiplicity, the corresponding Riesz projections converge in norm. Once \(\|\mathcal P_N^\lambda-\mathcal P^\lambda\|<1\), the two projections have the same finite rank. This proves multiplicity stability. Every nonzero isolated eigenvalue of a compact operator has finite algebraic multiplicity.
\end{proof}

\subsection*{Proof of Proposition~\ref{thm:a_posteriori_eigenvalue_bound}}
\begin{proof}
Set
\[
    \mathcal M^\lambda\coloneqq\mathcal P^\lambda\mathcal F,
    \qquad
    \mathcal M_N^\lambda\coloneqq\mathcal P_N^\lambda\mathcal F.
\]
The proof of Corollary~\ref{cor:spectral_projection} gives \(\mathcal P_N^\lambda\to\mathcal P^\lambda\) strongly. Since \(\mathcal M^\lambda\) is finite dimensional, the convergence is uniform on \(\mathcal M^\lambda\). Hence
\[
    S_N\coloneqq\mathcal P_N^\lambda|_{\mathcal M^\lambda}\colon
    \mathcal M^\lambda\to\mathcal M_N^\lambda
\]
is injective for all sufficiently large \(N\). Multiplicity stability gives equality of the dimensions, so \(S_N\) is an isomorphism. Define
\[
    J_N\coloneqq\mathcal P^\lambda|_{\mathcal M_N^\lambda}\colon
    \mathcal M_N^\lambda\to\mathcal M^\lambda.
\]
We have \(S_N\to I_{\mathcal M^\lambda}\) and
\[
    J_NS_N
    =\mathcal P^\lambda\mathcal P_N^\lambda|_{\mathcal M^\lambda}
    \longrightarrow I_{\mathcal M^\lambda}
\]
in norm. Thus \(S_N\) and \((J_NS_N)^{-1}\) are uniformly bounded. Since
\[
    J_N^{-1}=S_N(J_NS_N)^{-1},
\]
\(J_N\) is an isomorphism for all sufficiently large \(N\), with
\[
    \sup_{N\ge N_0}\|J_N^{-1}\|<\infty.
\]

Let
\[
    L\coloneqq\mathcal K|_{\mathcal M^\lambda},
    \qquad
    G_N\coloneqq\mathcal K_N|_{\mathcal M_N^\lambda}.
\]
Both spectral subspaces are invariant under their respective operators, and \(L\) has the single eigenvalue \(\lambda\), counted with algebraic multiplicity \(m_\lambda\). For \(x\in\mathcal M_N^\lambda\), commutation of the Riesz projections with their operators gives
\[
    (LJ_N-J_NG_N)x
    =\mathcal P^\lambda(\mathcal K-\mathcal K_N)x.
\]
Therefore,
\[
    \|LJ_N-J_NG_N\|
    \le\|\mathcal P^\lambda\|\,
       \|(\mathcal K-\mathcal K_N)\mathcal P_N^\lambda\|
    =\|\mathcal P^\lambda\|\,\eta_N.
\]

The eigenvalues of \(G_N\) are exactly the eigenvalues of \(\mathcal K_N\) inside \(\Gamma_\lambda\), counted with algebraic multiplicity. Hence
\[
    \operatorname{tr}(L)=m_\lambda\lambda,
    \qquad
    \operatorname{tr}(J_NG_NJ_N^{-1})
    =\operatorname{tr}(G_N)=m_\lambda\bar\lambda_N.
\]
Using \(|\operatorname{tr}(B)|\le m_\lambda\|B\|\) on the \(m_\lambda\)-dimensional space \(\mathcal M^\lambda\), we obtain
\[
\begin{aligned}
|\lambda-\bar\lambda_N|
&\le\|L-J_NG_NJ_N^{-1}\|\\
&\le\|\mathcal P^\lambda\|\,\|J_N^{-1}\|\,\eta_N.
\end{aligned}
\]
The uniform bound on \(\|J_N^{-1}\|\) proves Eq.~\eqref{eq:chatelin_residual_bound} with
\[
    C=\|\mathcal P^\lambda\|
      \sup_{N\ge N_0}\|J_N^{-1}\|.
\]
\end{proof}

\section{Additional pendulum diagnostics}\label{app:pendulum_additional}
This appendix gives supplementary pendulum calculations based on the coefficient-space matrices and trajectory data used in Figures~\ref{fig:pendulum_1} and~\ref{fig:pendulum_2}.

\begin{figure}
    \centering
    \includegraphics[width=0.95\linewidth]{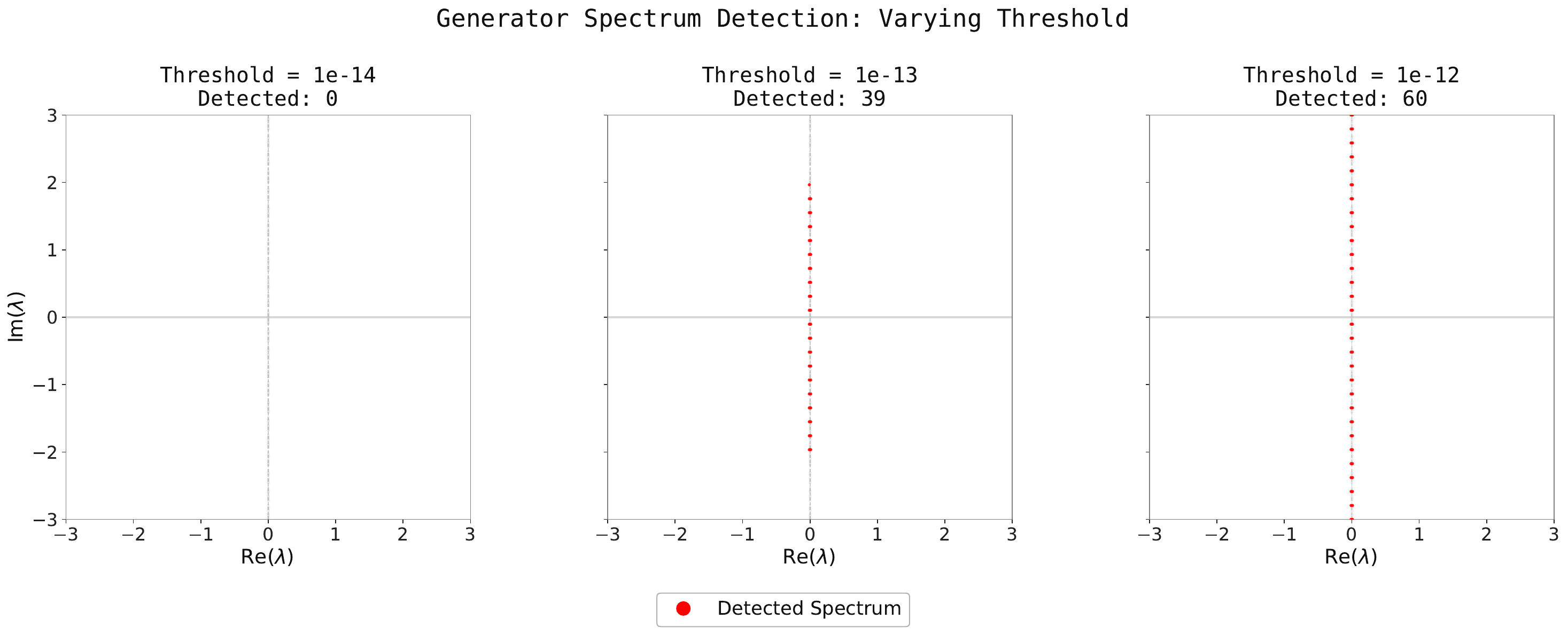}
    \caption{Coefficient-space threshold sensitivity of the forward-difference generator matrix on the prescribed lines \(\operatorname{Re} z=\pm0.01\), with \(\varepsilon\in\{10^{-14},10^{-13},10^{-12}\}\). The lines are fixed in advance and therefore do not independently localize the generator spectrum.}
    \label{fig:pendulum_3}
\end{figure}

\begin{figure}
    \centering
    \includegraphics[width=0.55\linewidth]{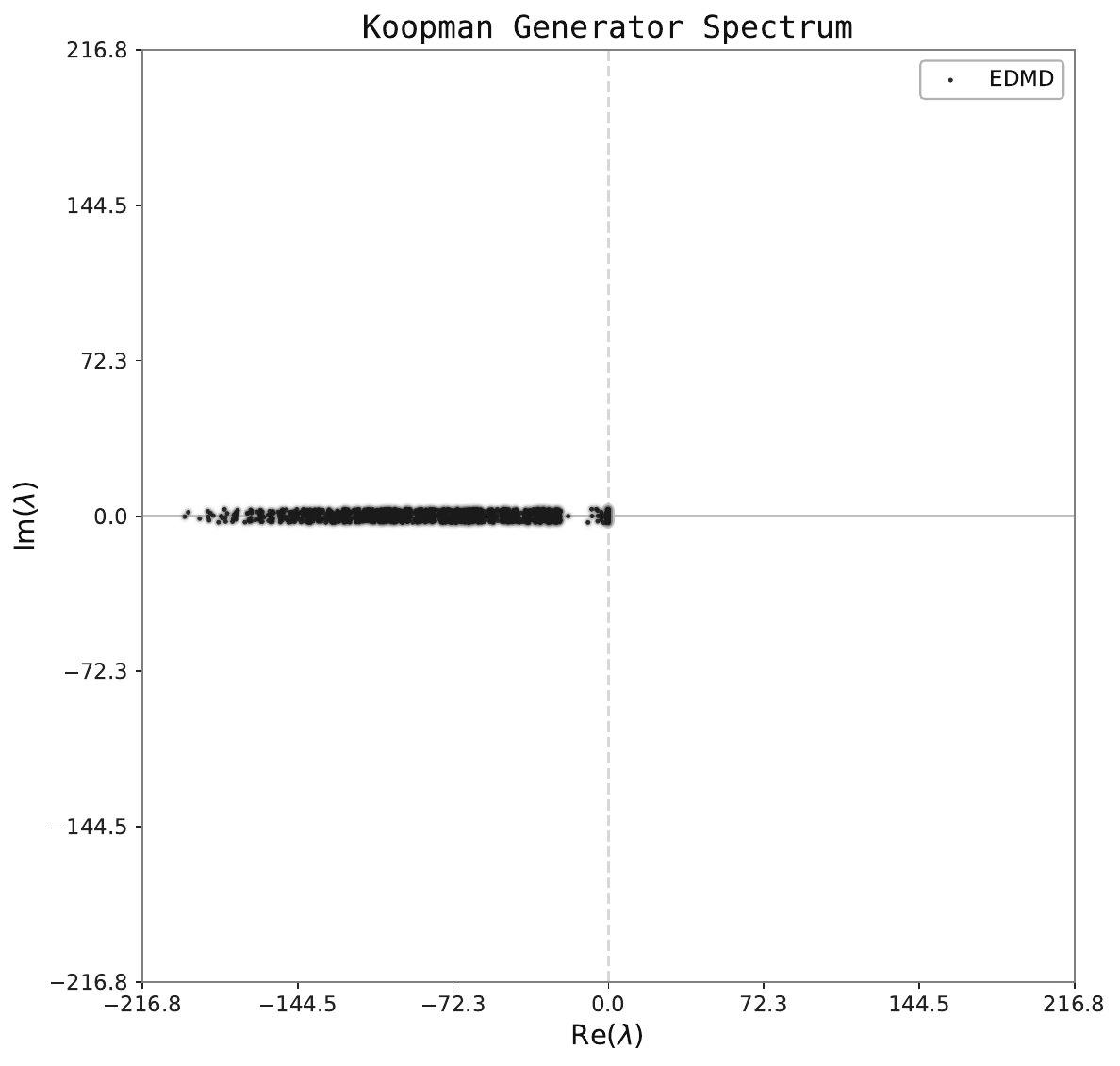}
    \caption{Eigenvalue-wise logarithmic map \(\lambda_j\approx\Delta t^{-1}\log\mu_j\) applied to the discrete EDMD eigenvalues of the Euler pendulum map, using the branch selected in the implementation. The result is branch dependent, subject to aliasing, and is not obtained from a matrix logarithm.}
    \label{fig:pendulum_4}
\end{figure}

\begin{figure}
    \centering
    \includegraphics[width=0.75\linewidth]{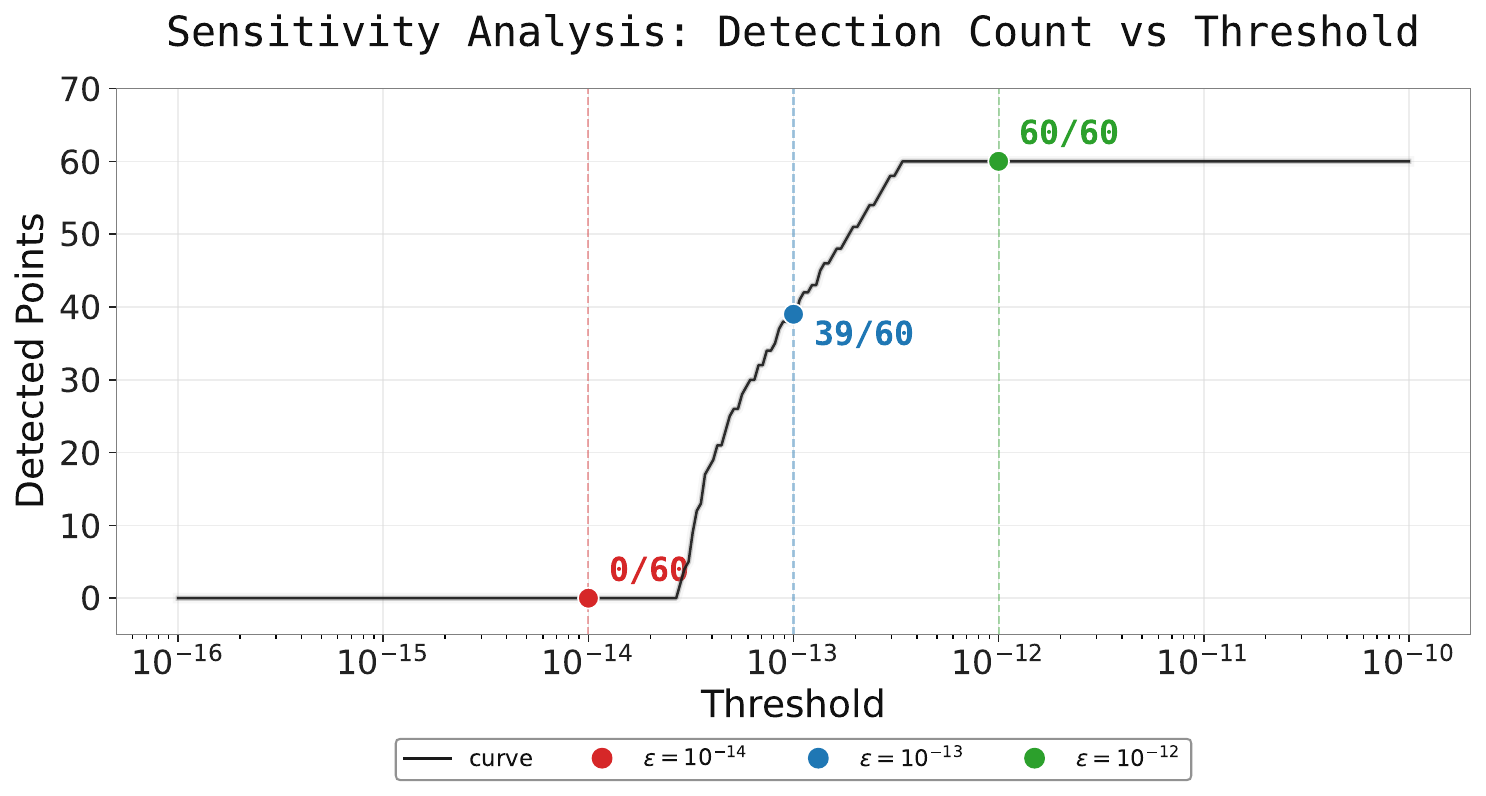}
    \caption{Number of sampled parameters satisfying the coefficient-space singular-value threshold as \(\varepsilon\) varies. Over the displayed range, the count changes from zero to all sampled points.}
    \label{fig:pendulum_5}
\end{figure}

Figure~\ref{fig:pendulum_3} evaluates the forward-difference generator matrix on the prescribed lines \(\operatorname{Re} z=\pm0.01\). Remark~\ref{rm:fd_accuracy} describes the associated temporal-discretization bias. Figure~\ref{fig:pendulum_4} applies the branch-dependent map \(\Delta t^{-1}\log\mu_j\) separately to the discrete EDMD eigenvalues. This is an eigenvalue-wise transformation rather than a matrix logarithm. Figure~\ref{fig:pendulum_5} records the number of selected spectral points as the coefficient-space threshold varies. The thresholds are close to double-precision resolution, so the calculation measures numerical sensitivity rather than a scale-independent spectral gap.

\section*{Acknowledgments}
Y.X. and I.I. acknowledge support from JST CREST Grant No.~JPMJCR24Q1, including the AIP Challenge Program. I.S. acknowledges support from JSPS KAKENHI Grant No.~JP24K20864.

\newpage
\bibliography{references}
\bibliographystyle{plain}

\end{document}